\documentclass[preprint, review,12pt,authoryear]{elsarticle}
\usepackage{graphicx}
\usepackage{fullpage,multirow, amsmath,amssymb}
\usepackage{array}
\newcolumntype{C}[1]{>{\centering\let\newline\\\arraybackslash\hspace{0pt}}m{#1}}
\newcolumntype{L}[1]{>{\raggedright\let\newline\\\arraybackslash\hspace{0pt}}m{#1}}
\newcolumntype{R}[1]{>{\raggedleft\let\newline\\\arraybackslash\hspace{0pt}}m{#1}}
\newcommand{\xvec}{\mathbf{x}}
\newcommand{\vvec}{\mathbf{v}}
\newcommand{\rvec}{\mathbf{r}}
\newcommand{\gvec}{\mathbf{g}}
\newcommand{\pvec}{\mathbf{p}}
\newcommand{\dvec}{\mathbf{d}}
\newcommand{\yvec}{\mathbf{y}}
\newcommand{\evec}{\mathbf{e}}
\newcommand{\fref}[1]{Fig.~\ref{#1}}
\makeatletter
\def\Ginclude@eps#1{%
 \message{<#1>}%
  \bgroup
  \def\@tempa{!}%
  \dimen@\Gin@req@width
  \dimen@ii.1bp%
  \divide\dimen@\dimen@ii
  \@tempdima\Gin@req@height
  \divide\@tempdima\dimen@ii
    \includegraphics{#1}%
  \egroup}
\makeatother

\journal{Journal of Petroleum Science and Engineering}

\begin{document}
\begin{frontmatter}

\title{Joint optimization of well placement and control for nonconventional well types}

\author[osu]{T. D. Humphries\corref{cor1}}
\ead{humphrit@science.oregonstate.edu}
\author[mun]{R. D. Haynes}
\ead{rhaynes@mun.ca}

\address[osu]{Department of Mathematics, Oregon State University, Corvallis, OR, USA}

\address[mun]{Department of Mathematics and Statistics, Memorial University of Newfoundland, St. John's, NL, Canada}
\cortext[cor1]{Corresponding author. +1-541-737-5165}

%\maketitle

\begin{abstract}

Optimal well placement and optimal well control are two important areas of study in oilfield development. Although the two problems differ in several respects, both are important considerations in optimizing total oilfield production, and so recent work in the field has considered the problem of addressing both problems jointly. Two general approaches to addressing the joint problem are a simultaneous approach, where all parameters are optimized at the same time, or a sequential approach, where a distinction between placement and control parameters is maintained by separating the optimization problem into two (or more) stages, some of which consider only a subset of the total number of variables. This latter approach divides the problem into smaller ones which are easier to solve, but may not explore search space as fully as a simultaneous approach.

In this paper we combine a stochastic global algorithm (Particle Swarm Optimization) and a local search (Mesh Adaptive Direct Search) to compare several simultaneous and sequential approaches to the joint placement and control problem. In particular, we study how increasing the complexity of well models (requiring more variables to describe the well's location and path) affects the respective performances of the two approaches. The results of several experiments with synthetic reservoir models suggest that the sequential approaches are better able to deal with increasingly complex well parameterizations than the simultaneous approaches.
\end{abstract}

\begin{keyword}
Production Optimization \sep Waterflooding \sep Well Placement \sep Well Control \sep Mesh Adaptive Direct Search \sep Particle Swarm Optimization \sep Derivative-free Optimization \sep Sequential Optimization
\end{keyword}
\end{frontmatter}

\section{Introduction}

Maximizing oilfield production is a problem of considerable economic importance, which has received attention from both the industrial and academic communities. Given a computational model of the reservoir, one typically considers drilling and operating two types of well. Production wells (or producers) withdraw fluid from the reservoir, while injection wells (injectors) pump either water or gas into the reservoir, to drive oil towards producers once the incipient pressure in the reservoir has dissipated. The pumping rates of these wells can be either controlled directly, or indirectly by prescribing a bottom hole pressure (BHP) at each well. By holding injectors at a higher BHP than producers, fluid is driven from the former to the latter; the actual flow rate achieved depends on the size of the pressure gradient that is created, as well as the reservoir geology. 

In the context of optimization, the relevant decision variables may include the number and type of wells to drill, the order in which to drill them, their positions and orientations, and the control parameters (pumping rates or BHPs) for each well. Until recently, most work on production optimization has focused on optimizing only a subset of these variables at one time. The two most frequently studied problems are determining optimal well locations (well placement), and determining optimal control parameters for wells already in place. In both cases, the optimization problem is challenging for several reasons: evaluating the objective function requires computationally expensive reservoir simulations; gradient information about the objective function is often unavailable (or expensive to approximate); and the function is generally nonconvex and may include many local optima. Depending on the problem being considered, the number of variables may range from less than ten to more than one hundred.

A typical objective function used in production optimization problems is the net present value (NPV) of produced oil. Let $\xvec$ denote the vector of parameters being optimized (well positions, controls, etc.). A standard formula for NPV (used in this paper) is
\begin{align}
NPV(\xvec) = C(\xvec) + \int_0^T  \Biggl\{  &\sum_{n \in prod} \left[  c_o q_{n,o}^- (\xvec, t) - c_{w,disp} q_{n,w}^- (\xvec, t) \right] \notag \\
			&-   \sum_{n \in inj} c_{w,inj} q_{n,w}^+(\xvec, t) \Biggr\} (1+r)^{-t}\: dt, \label{E:NPV}
\end{align}
\noindent where $t \in [0, T]$ represents time over the production period of $T$ years. The parameters $c_o$, $c_{w,disp}$ and $c_{w,inj}$ represent the price per barrel of produced oil, disposal cost per barrel of produced water, and cost per barrel of injected water, respectively, and the annual interest rate is specified by $r$. The functions $q_{n,o}^-(\xvec, t)$ and $q_{n,w}^-(\xvec, t)$ are the production rates (barrels/day) of oil and water, respectively, at well $n$, while $q_{n,w}^+(\xvec, t)$ is the water injection rate at well $n$; $n$ must be a member of either the set of producers ($prod$) or injectors ($inj$).  The rates are determined from a reservoir simulation under the operating parameters described by $\xvec$. $C(\xvec)$ represents capital expenditures such as drilling costs, and may be omitted if these costs are constant with respect to the optimization parameters. Although the NPV formula prioritizes producing more oil early in the production period (for $r > 0$), other considerations exist as well. For example, once the water front from an injector arrives (breaks through) at a producer, the well is flooded and begins producing large amounts of water. In a reservoir with multiple water injectors and producers, one typically tries to delay breakthrough of a water front at a producer until fronts from other nearby injectors have also arrived, thus increasing the amount of oil that the well produces prior to being flooded.

%Along these lines, one problem that has received recent attention is that of jointly optimizing well placement and well control. 
%
%Existing literature on the joint optimization problem suggests that a fully simultaneous approach provides better results than a sequential one. Most of this work has focused only on placing vertical wells, whose positions depend only on two parameters. In many applications, however, it is of interest to consider wells with more complex parameterizations, such as horizontal or arbitrarily oriented wells.

%\section{Background}\label{S:background}

Well placement problems involve optimizing over parameters corresponding to the positions and orientations of the injection and production wells. Typically these problems assume a simple reactive control scheme, where injectors and producers are held at fixed BHPs during the entire production period, and producers are eventually shut in once the ratio of water to oil produced exceeds a profitable threshold.  Optimal well placement depends in large part on the permeability field of the reservoir, which tends to be highly heterogeneous; as a result, the objective function surface in well placement problems is fairly rough, and may include many local optima~\citep{CMD11,OD10}. Approaches based strictly on following the gradient of the objective function will typically converge to the closest local optimum without fully exploring the solution space; for this reason, studies on optimal well placement have tended to focus on stochastic, heuristic approaches aimed at exploring the solution space globally. Historically, genetic algorithms (GAs) have received the widest use~\citep{BH97,YDA03,GH04,ADOA06,ESMAS09}.  More recently, other stochastic approaches such as simultaneous perturbation stochastic approximation (SPSA)~\citep{BKWSS06}, covariance matrix adaptation~\citep{BDA12}, particle swarm optimization (PSO)~\citep{OD10,OD11, DJLT14} and differential evolution (DE)~\citep{NNR13} have also been applied successfully to the problem. Gradient-based optimization approaches have been successfully adapted to the well placement problem in some cases as well~\citep{ZHVBJ08,ZLRY10}.

In well control optimization problems, the decision variables usually consist of pumping rates or BHPs. One often assumes that these parameters are held constant during several time intervals of fixed length, between which they can altered (e.g. every two years). Thus the total number of control parameters to be determined is the product of the number of wells and the number of time intervals. In comparison with the well placement problem, the NPV function varies much more smoothly when control parameters are altered. Thus, optimization approaches which are based on approximating the gradient using the adjoint method have been the most popular approach to this problem (see~\citet{J11} for a recent review). One drawback of adjoint approaches is that the reservoir simulator must provide adjoint information, which is not always available. For this reason, we will consider only black box optimization algorithms, which are based entirely on simulator output and do not require an explicit gradient calculation.  Examples of black-box algorithms that have been applied to the well control problem include SPSA~\citep{WLR09, DR13}, methods based on quadratic models and trust regions~\citep{HAM13, ZCODLR13}, genetic algorithms~\citep{YPKS13}, augmented Lagrangian methods~\citep{CWLR10}, and generalized pattern search (GPS) approaches~\citep{CID10,CMD11, ANDK14}.

In general, the optimal positioning of wells depends to some extent on the control scheme that is employed. Thus, recent work has considered jointly optimizing well placement and well control, with the expectation that one can find better solutions than by a purely placement-based approach. One proposed method~\citep{FLR10, FR13} places a large number of injectors and producers in the reservoir initially, and then uses an adjoint method to determine optimal well controls. Wells with low flowrates are removed from the simulation, thereby determining optimal positions as well. A second approach~\citep{BCDFK12} uses a combination of GPS and adjoint methods in a nested optimization procedure. The outer iteration consists of using GPS to determine optimal well positions, where the objective function involves performing an inner optimization using the adjoint method to determine the best control strategy. The SPSA algorithm has also been applied to problems involving both well placement and well control in~\citet{LJ12,LJM13}. Other papers have considered more specialized procedures such as determining optimal well positions using a heuristic screening process and then optimizing control~\citep{XZCL13}, and using a MINLP approach that is fully integrated with the reservoir simulator~\citep{TKTBA13}. Finally, the use of PSO and pattern search in tandem has been investigated in~\citet{IDE13, IED14}, and by the authors of this paper in~\citet{HHJ13}; the first two papers used a variant of pattern search known as Mesh Adaptive Direct Search (MADS).

The wide variety of algorithms that have been applied to different production optimization problems is reflective of the characteristics of the objective function in these different problems. Local searching methods based on adjoint gradients or quadratic models are able to exploit smooth features of the objective function, and converge rapidly to a local optimum. Stochastic approaches such as GAs, PSO and DE converge more slowly, but their focus on global searching is well-suited to highly multimodal objective functions. Methods such as GPS and MADS perform some degree of global and local searching, without requiring any gradient approximations. Thus it is natural that approaches to the joint problem, where the objective function has smooth behaviour with respect to some variables and rough behaviour with respect to others, tend to combine several of these algorithms.

An important question in the context of the joint problem is whether it is preferable to optimize over all parameters simultaneously, or to use a sequential approach that divides the problem into smaller subproblems. Optimizing over all variables simultaneously ensures that the best possible solution exists somewhere in the search space. This search space may be very large, however, and given the highly nonconvex nature of the optimization problem, it may be difficult to find the global optimum. One can instead consider a sequential procedure; for example, determining optimal well positions first while assuming some simple control scheme, and then using the best solution found as an initial guess for the full problem. The first subproblem involves fewer variables than the full problem and should therefore be easier to solve; however, it is possible that the configuration of wells found by solving the first subproblem (which will depend on the assumed control scheme) is not close to the global optimum of the full problem. Thus the second step of the optimization may not be able to find the global optimum.

Several of the aforementioned papers addressing the joint problem~\citep{LJM13, BCDFK12,IDE13} have found that the solutions found by optimizing over all variables simultaneously are superior to those that can be found by sequential approaches. In~\citet{HHJ13}, however, we found that a sequential approach (denoted as the ``decoupled'' approach in that paper) was competitive with the studied simultaneous approach, and was even preferable in some test cases. Thus, we wish to further investigate the effectiveness of these two approaches under different experimental conditions. In particular, most papers on the topic have considered only drilling vertical wells, which are parameterized simply by their $(x,y)$ co-ordinates. In many applications, however, it is of interest to consider drilling wells with more complex parameterizations, such as horizontal or arbitrarily oriented (deviated) wells~\citep{YDA03,OD10}. In these cases it may be more difficult to determine the optimal positioning of wells, meaning that focusing on well placement before considering the joint placement and control problem could be advantageous.

The paper is organized as follows. In Section~\ref{S:methods} we describe the primary optimization algorithms that are used (Mesh Adaptive Direct Search and Particle Swarm Optimization, as well as a hybridization of the two), and how they are applied to the problem. The optimization framework is implemented using the NOMAD black-box software~\citep{D11, AACDD}. In Section~\ref{S:experiments} we describe our numerical experiments. Results are presented in Section~\ref{S:results} and discussed in Section~\ref{S:discussion}. We finish with concluding remarks in Section~\ref{S:concl}. Overall the results of our experiments suggest that our sequential approach is more robust than a fully simultaneous approach in cases where well parameterizations are more complex.

\section{Methodology}\label{S:methods}
In this section we briefly describe the optimization approaches used in our experiments. These algorithms are well-suited to the production optimization problem because they require no gradient information and are easily parallelizable, which helps offset the high cost of function evaluations. We also briefly discuss handling of nonlinear constraints which arise in production optimization problems.

\subsection{Mesh Adaptive Direct Search} \label{SS:MADS}

Mesh Adaptive Direct Search (MADS)~\citep{AD06} is a pattern search algorithm which extends the earlier class of Generalized Pattern Search (GPS) algorithms \citep[e.g.][]{LT99,AD03}. A GPS algorithm is an iterative algorithm consisting of a series of {\em search} and {\em poll} steps. Let $N$ denote the number of optimization variables. At every iteration $k$, a discrete mesh, centred at the current incumbent point $\xvec^{(k)} \in \mathbb{R}^N$, is defined by:
\begin{equation}
M^{(k)} =  \left\{ \xvec^{(k)} + \Delta^m_k D \mathbf{z} \: : \: \mathbf{z} \in \mathbb{N}^{n_D} \right\}, \notag
\end{equation}
where $\Delta^m_k$ controls the resolution of the mesh at iteration $k$, $D$ is a matrix whose columns consist of the polling directions, $\mathbb{N}$ is the set of natural numbers, and $n_D$ is the number of polling directions. The polling directions must form a positive spanning set in solution space; i.e., one must be able to specify any point in solution space by adding together only positive scalar multiples of these directions. A common choice of directions is
\begin{align}
\mathcal{D} &= \left\{ \evec_1, \evec_2, \dots, \evec_N, -\evec_1, -\evec_2, \dots, -\evec_N \right\}, \notag
\end{align}
where the $\evec_n$ are the canonical basis vectors $(1,0,0, \dots, 0)^T$, $(0,1,0, \dots, 0)^T$, etc.  Here $\mathcal{D}$ refers to the set of polling directions, which form the columns of the matrix $D$.

The {\em search} step is a generic optimization step in which the objective function is evaluated at some finite number of points on $M^{(k)}$. The number of points and their co-ordinates can be chosen according to any desired strategy; the only requirements are that the points lie on the mesh and can be computed in a finite amount of time. If the best point found overall is an improvement on $\xvec^{(k)}$, then it becomes the new incumbent; otherwise, the algorithm proceeds to the {\em poll} step.  The poll step consists of evaluating the objective function at all the points that are immediate neighbours of the incumbent point on the mesh $M^{(k)}$. These points are given by
\begin{equation}
\left\{ \yvec_j^{(k)} \right\} = \left\{ \xvec^{(k)} + \Delta^m_k \dvec_j | \: \forall \: \dvec_j \in \mathcal{D} \right\}.\label{E:MADS}
\end{equation}

If the best point found overall by polling is an improvement, then it becomes the new incumbent. $\Delta^m_k$ may then be increased for the next iteration. If the poll step is unsuccessful, then $\Delta^m_k$ is reduced and another iteration begins, using the same incumbent point as before. The algorithm is considered to have converged once $\Delta^m_k$ is reduced beyond some minimum threshold, which indicates that the current point is at least close to a local optimum. In fact, provided that the objective function is continuously differentiable, GPS is guaranteed to converge to a local optimum, at least to mesh precision~\citep{LT99}.

The key distinction between GPS and MADS is that MADS introduces a second parameter $\Delta^p_k$, which controls the polling size and is required to be larger than the mesh size $\Delta^m_k$, for all $k$. Thus the size of the underlying mesh $M^{(k)}$ shrinks more quickly than the size of the polling stencil as the algorithm proceeds. This allows MADS to select the set of polling directions, $\dvec_j \in \mathcal{D}$, on a finer grid than GPS does, and thereby generate an asymptotically dense set of polling directions as $\Delta^m_k$ and $\Delta^p_k$ decrease~\citep{AD06}. Unlike in GPS, $\mathcal{D}$ can be changed with every MADS iteration, and thus MADS is not limited to a finite number of polling directions as the mesh size decreases. MADS was found to provide superior results for several test problems in~\citet{AD06} as a result. There are several strategies for generating $\mathcal{D}$; we use OrthoMADS~\citep{AADD09}, which deterministically generates a set of $2N$ orthogonal polling directions at every iteration.

\subsection{Particle Swarm Optimization} \label{SS:PSO}

Particle swarm optimization (PSO)~\citep{KE95,C06} is an iterative heuristic algorithm based on collective intelligence. The swarm consists of some number of particles (typically 20--50), with the $i$th particle being characterized by its position $\xvec_i^{(k)}$ and velocity $\vvec_i^{(k)}$ at every iteration $k$. Both $\xvec_i^{(k)}$ and $\vvec_i^{(k)}$ are vectors of size $N$. The position corresponds to a vector in search space, with an associated objective function value $f(\xvec)$. Every particle retains memory of the best position it has found so far, denoted $\pvec_i^{(k)}$, and communicates with other particles in some neighbourhood to determine the best position found among all of them, denoted $\gvec_i^{(k)}$. Following initialization, the algorithm proceeds as follows:

\begin{align}
\xvec_i^{(k+1)} &= \xvec_i^{(k)} + \vvec_i^{(k+1)}, \label{E:PSO1} \\
\vvec_i^{(k+1)} &= \iota \vvec_i^{(k)} + \mu \rvec_1^{(k)} \otimes \left(\pvec_i^{(k)} - \xvec_i^{(k)} \right) +  \nu \rvec_2^{(k)} \otimes \left(\gvec_i^{(k)} - \xvec_i^{(k)} \right).\label{E:PSO2} 
\end{align}

The velocity update (\ref{E:PSO2}) combines three terms; a tendency to continue moving in the direction given by the particle's current velocity, $\vvec_i^{(k)}$, a tendency to be drawn towards the best position the particle has found so far, $\pvec_i^{(k)}$, and a tendency to be drawn towards the best position found by all particles in its neighbourhood, $\gvec_i^{(k)}$. The constants $\iota$, $\mu$ and $\nu$ are weighting parameters. The $N$-vectors $\rvec_1^{(k)}$ and $\rvec_2^{(k)}$ are randomly generated from the uniform distribution on $(0,1)$ at every iteration, with $\otimes$ denoting componentwise multiplication. This multiplication adds a stochastic component to the algorithm and helps avoid early convergence to a local minimum. PSO is typically run for some fixed number of iterations or until some convergence criterion has been satisfied, e.g., until the particle velocities are close to zero. In general, there is no guarantee that PSO converges to a global or even a local optimum; however, in practice it has proved to be effective for a wide variety of optimization problems~\citep{P08}. 

Following~\citet{C06}, our implementation of PSO uses standard weighting parameters of $\iota = 0.721$, and $\mu = \nu = 1.193$ and a population size of 50 particles. We use a global best neighbourhood topology, meaning that  every particle communicates with every other particle in the swarm, and thus $\gvec_i^{(k)}$ can be replaced by a single vector $\gvec^{(k)}$, representing the best solution found so far. Although it has been observed that this choice of neighbourhood may cause PSO to converge to a local optimum before thoroughly exploring search space~\citep{C06}, we found that other choices of neighbourhood caused the algorithm to converge too slowly for our purposes. Our implementation of PSO also includes some heuristics for handling nonlinear constraints, which are discussed in Section~\ref{SS:constr}. 

\subsection{MADS-PSO} \label{SS:MADS-PSO}
Generally speaking, stochastic approaches such as PSO are well-suited to problems where the objective function surface is rough, while MADS is better-suited to problems where the objective function varies more smoothly. Since the joint well placement and control optimization problem features both types of behaviour, one expects that combining the two algorithms could be beneficial. This has been shown to be true in previous studies~\citep{IDE13, HHJ13}. The two algorithms can be readily hybridized by implementing PSO as the {\em search} step within the GPS framework, as was first proposed in~\citet{VV07,VV09}. During this search step, the particles move according to the update formulas~(\ref{E:PSO1}) and  (\ref{E:PSO2}), and are then projected onto the mesh. If the search step fails to improve on the incumbent solution, polling takes place around the current best position found, as per (\ref{E:MADS}).  If the poll step finds a better solution, the current best position is updated and a new iteration begins; otherwise, the polling stencil size is reduced.

\subsection{Sequential approaches} \label{SS:SEQ}
The MADS, PSO, and MADS-PSO approaches described previously can all be applied to the joint problem by optimizing over all control and placement parameters simultaneously. In contrast, our sequential approach to the joint optimization problem splits it into a two-step procedure:
\begin{enumerate}
\item Optimize well positions using PSO, assuming some simple control scheme for injectors and producers.
\item Use the best position found during Step 1 as the initial incumbent solution for optimizing over placement and control simultaneously using MADS.
\end{enumerate}
We note that our approach differs from the sequential approaches studied elsewhere in the literature \citep[e.g.][]{BCDFK12,IDE13} where well positions are held fixed during the second step, and only well controls are optimized. There is no drawback to allowing well positions to be further optimized during the second step, aside from some increased computational cost due to the higher dimensionality of the problem. In our experiments we have found that the increased cost is justified, as it is possible to find significantly better solutions during the second step of optimization if one allows the positional parameters to vary in addition to the control parameters. In essence, then, the second step of our sequential approach consists of optimizing over all variables simultaneously using MADS, starting from a good initial guess which has been provided by solving a smaller well placement optimization problem in the first step. 

The choice of control scheme for Step 1 has some bearing on the performance of this approach, as previously shown in~\citet{HHJ13}. We will therefore study two variants:
\begin{itemize}
\item Sequential-I: Injectors are held at the maximum allowable BHP during the entire production period, while producers are held at the minimum allowable BHP. This control scheme is a common choice for well placement problems, and generates the highest flow rates possible for a given configuration of injectors and producers.
\item Sequential-II: Injectors are held at 10-20\% below the maximum BHP, while producers are held 10-20\% above the minimum BHP. The idea is that this choice of controls might better approximate the average BHP of a well over the whole production period when the controls can be altered. Thus, the well positions found during the first step of optimization may be closer to the optimal configuration for the full problem, which could improve the odds of finding a globally optimal solution in the second step.
\end{itemize}

\subsection{Constraint handling} \label{SS:constr}
Well placement and control optimization problems always include bound constraints, since wells must lie within the boundaries of the reservoir, and are subject to operating constraints that restrict the values of the control parameters. General constraints are often present as well. A common example is a limit on fluid injection and production rates for wells controlled by BHP. In this case the well flow rates, which may have minimum and maximum values depending on operating constraints, have a complicated nonlinear dependence on the reservoir geology and positions of the various wells, in addition to the prescribed BHPs. A violation of one of these constraints can be quantified using a constraint violation function; for instance, if there are maximum flowrates for injectors and producers (denoted $q^{max}_{inj}$ and $q^{max}_{prod}$, respectively), then one can use
\begin{align}
h(\xvec) =  \sum_{n \in prod} &\left \{ \int_0^T \max \left( q_{n,o}^- (\xvec, t) + q_{n,w}^- (\xvec, t) -q^{max}_{prod}, 0 \right) \: dt  \right\} \notag \\
					&+ \sum_{n \in inj} \left\{ \int_0^T \max \left( q_{n,w}^+(\xvec, t)-q^{max}_{inj}, 0\right)\: dt \right \},\label{E:constr}
\end{align}
Thus $h=0$ if a solution completely satisfies the constraints, and is positive otherwise; a point that violates the constraints is {\em infeasible}.

Both the PSO and MADS components of our optimization approaches must therefore include some method of dealing with bound and general constraints. For MADS, we use default options that are provided with NOMAD. Points that violate bound constraints are projected back onto the boundary of search space, while general constraint violations are handled using the progressive barrier approach~\citep{ADD10}. This approach permits infeasible points to be considered as incumbents, provided that the constraint violation is below some threshold $h^{max}$, which decreases as the iteration proceeds.

For PSO, particles that travel outside the boundaries of search space are also brought back to the boundary, and have their velocities altered to avoid traveling outside the boundary again at the next iteration~\citep{C06}. To handle general constraints, we establish the following ranking system for determining each particle's personal best position as well as the global best:
\begin{enumerate}
\item Between any two infeasible solution, the one with the smaller $h$ value is best,
\item Any feasible solution is better than an infeasible solution, and
\item Between two feasible solutions, the one with the better objective function value is best.
\end{enumerate}
This approach, which was proposed in~\citet{CEC08}, has the advantage of being simple to implement and not requiring any parameter tuning (unlike penalty function based approaches), while also allowing particles to be initialized to infeasible positions.

Other general constraints that arise in the production optimization problem include the following:
\begin{enumerate}
\item Two wells must not intersect the same grid cell,
\item The well path should not intersect with inactive grid cells, and
\item Wells should be separated from one another by some minimum distance.
\end{enumerate}
All of these conditions can be checked prior to running a reservoir simulation. Points which violate the first two conditions are usually treated as invalid input and are not considered as potential solutions. The well distance constraint can be treated similarly, although it may be advantageous to instead assign it a numerical penalty based on the extent of the constraint violation, as in Equation (\ref{E:constr}). The latter approach would allow searching around promising infeasible solutions which only violate the distance constraint to a small extent. We did not include a well distance constraint in any of our experiments, as the optimal solutions found by our algorithm generally featured adequately-spaced wells, even without any explicit constraint.

Finally, some trial points may cause a reservoir simulation job to fail; for instance, if the numerical solver does not converge for a given set of well positions and controls. The NOMAD optimization software which was used to implement our optimization approach is capable of catching failed jobs and ensuring that the optimization run itself does not hang or crash. A set of parameters that causes the simulator to fail is essentially treated as invalid input, and thus not considered as a potential solution.

\section{Experiments} \label{S:experiments}
\subsection{Optimization framework}
The optimization approaches studied in this paper have been implemented using NOMAD~\citep{D11, AACDD}. This black-box optimization package includes an implementation of MADS, and also provides users with the ability to implement their own custom approaches to be used as a Search component. We have therefore written an implementation of PSO that can be used in conjunction with MADS (for the MADS-PSO algorithm), and also run as a straightforward PSO algorithm, for the purposes of comparison. The IMEX Advanced Oil/Gas Reservoir Simulator~\citep{CMG} has been used as the reservoir simulator, with interfacing between NOMAD and IMEX handled by custom code written in Python. We have also taken advantage of NOMAD's MPI-based implementation to evaluate the objective function at multiple trial points in parallel, by running the software on a computational cluster consisting of 48 nodes.

\subsection{Experiment 1}\label{SS:exp1}
\begin{figure}
\includegraphics[width=0.75\linewidth]{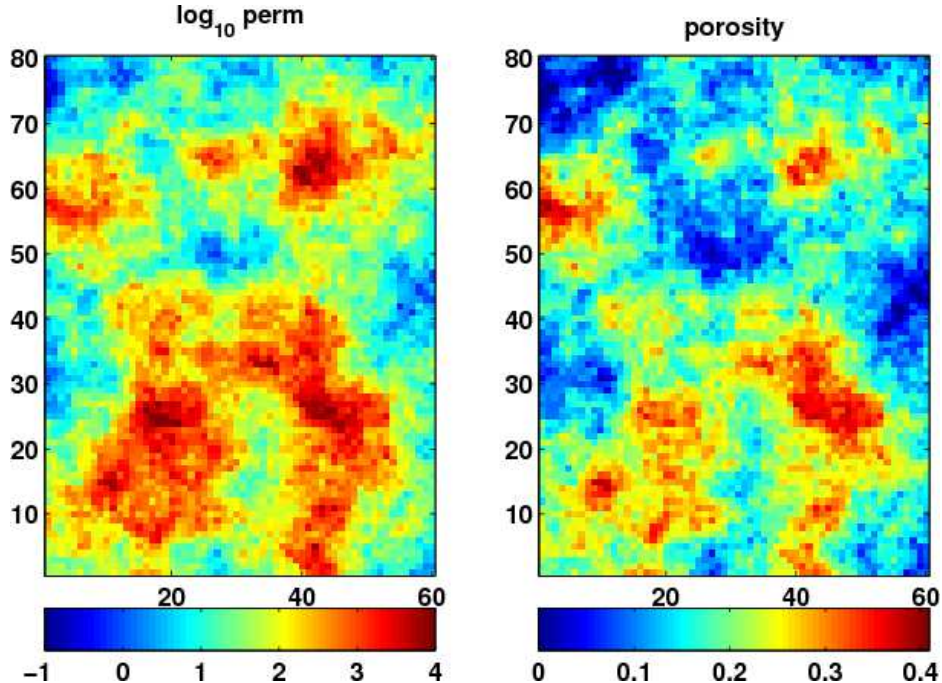}
\caption{Permeability (left) and porosity (right) fields used for Experiments 1 and 2. Permeability is shown in millidarcy (mD), on a logarithmic scale.}\label{F:perm}
\end{figure}

Our first test case consists of placing six vertical wells (two injectors and four producers) in a single-layer synthetic reservoir consisting of 60 $\times $ 80 cells, assuming a two-phase oil/water fluid component model. The initial saturation is 80\% oil to 20\% water. The reservoir permeability and porosity fields are shown in \fref{F:perm} and are taken from the~\citet{SPE10}. The production period consists of ten years with a control interval of two years. Thus each well has two positional parameters (its $(x,y)$ co-ordinates) and five control parameters, for a total of 42 variables. The $(x,y)$ co-ordinates correspond to grid indices, and each vertical well is assumed to be drilled in the centre of the cell. Even if one disregards the control parameters, there are ${4800 \choose 6} \approx 1.7 \times 10^{19}$ possible ways to position the six wells, so an exhaustive search is out of the question. The simulation parameters are summarized in Tables~\ref{T:exp_params} and \ref{T:fluid}, and the economic parameters used to compute NPV are given in Table~\ref{T:econ}. We considered two cases: Case 1A has no constraints on production, while Case 1B incorporates a nonlinear constraint by imposing a maximum flowrate of 1500 m$^3$/day on injectors and 750 m$^3$/day on producers. Equations~(\ref{E:NPV}) and (\ref{E:constr}) are used to compute the NPV and constraint violation $h$, respectively. We consider solutions with $h$ values slightly larger than zero to be feasible, to account for the fact that well flow rates may spike momentarily when bottom hole pressures are instantaneously changed from one value to another. We view this as an artifact of the simulation and not a true constraint violation.

\begin{table*}
\caption{Parameters for the three experiments considered.}\label{T:exp_params}
\begin{tabular}{lrrr}
			&Experiment 1			&Experiment 2 			&Experiment 3 \\
\hline \noalign{\smallskip}
Reservoir grid (cells)	&60 $\times$ 80 $\times$ 1	&60 $\times$ 80 $\times$ 1	&60 $\times$ 50 $\times$ 3 \\
Reservoir size (m)	&2400 $\times$ 3200 $\times$ 10		&2400 $\times$ 3200 $\times$ 10 	& 2400 $\times$ 2000 $\times$ 75 \\
Reservoir depth (m) &2000						&2000						&2000 \\
Oil-water contact (m) &2008						&2008						&2060 \\
Number of wells	&2 inj, 4 prod 			&2 inj, 2 prod 					&2 inj, 2 prod \\
Injector BHP bounds (bar)  &300 -- 450 &300 -- 450 &300 -- 450 \\
Producers BHP bounds (bar)  &125 -- 260  &125 -- 260 &125 -- 260 \\
Production period (years) 	     &10 		&10 	&15  \\
Control interval (years) 	     &2 			& 2 	& 3 \\
Positional parameters per well &2			&4		&6 \\
Control parameters per well & 5			&5		&5 \\
Total number of parameters &42 			&36 		&44
\end{tabular}
\end{table*}

\begin{table*}
\caption{Fluid properties used for all three experiments.} \label{T:fluid}
\begin{tabular}{lr}
Property			&Value \\
\hline
Water and oil density $\rho_{w}$, $\rho_{o}$		&1000 and 860 kg/m$^3$ \\
Water and oil viscosity $\mu_{w}$, $\mu_{o}$		&0.32 and 0.53 cp at 280 bars \\
Water and oil compressibility $c_w$, $c_o$		&$5\times 10^{-5}$ and $4.35\times10^{-5}$ bar$^{-1}$ \\
Relative permeability			&See curves below\\
\end{tabular}

\includegraphics[width=0.35\linewidth]{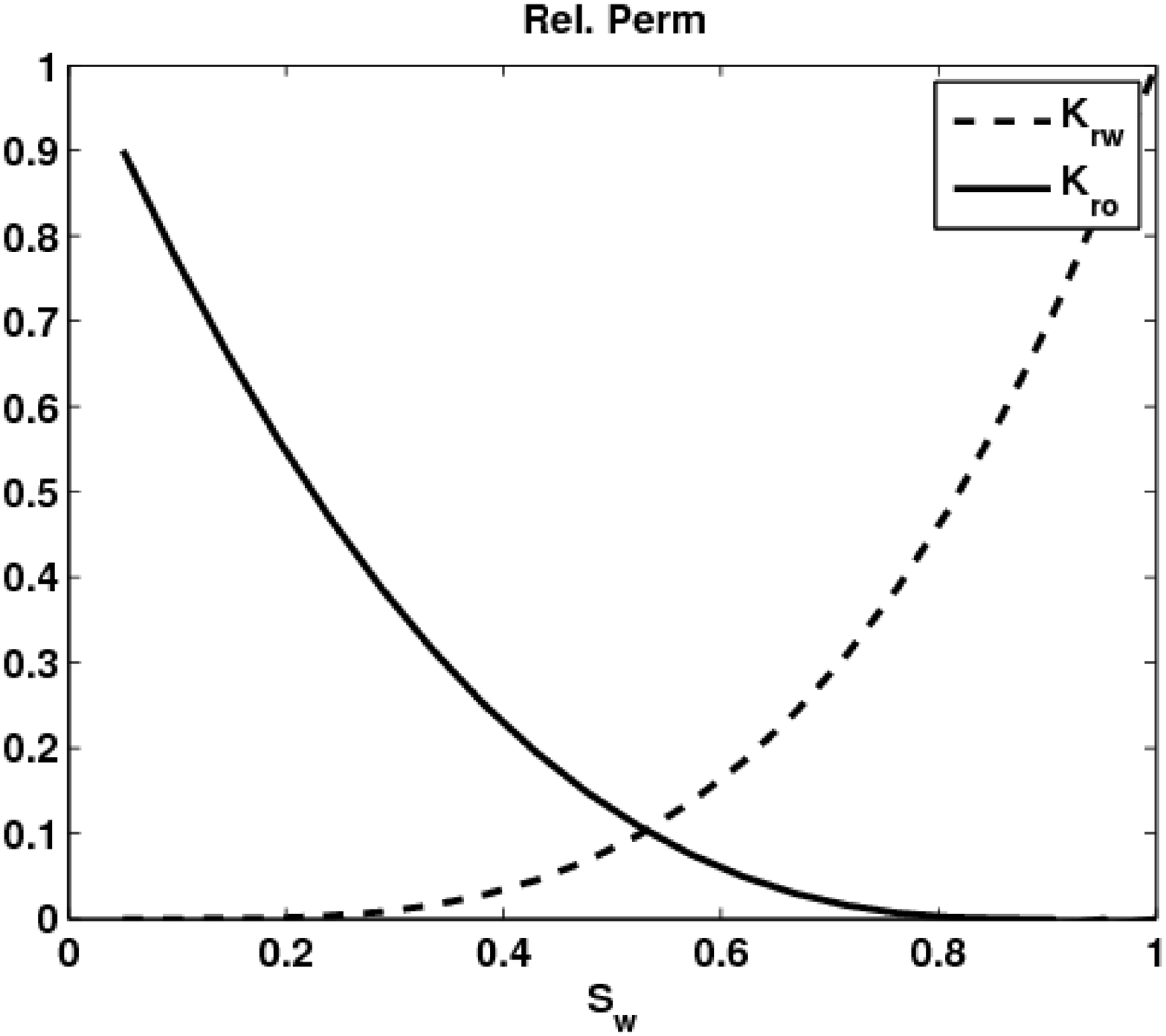}
\end{table*}

\begin{table}
\caption{Economic parameters used for all three experiments.}\label{T:econ}
\begin{tabular}{lr}
Parameter	&Values \\
\hline\noalign{\smallskip}
Value of produced oil, $c_o$	&\$80/barrel \\
Cost of injecting water, $c_{w,inj}$	&\$8/barrel \\
Cost of disposing of water, $c_{w,disp}$ &\$12/barrel \\
Interest rate, $r$			&10\% per annum \\
Base drilling cost	&\$25M per well\\
Drilling cost per unit length		&\multirow{2}{*}{\$50K per metre}\\
 (Experiments 2 and 3 only)
\end{tabular}
\end{table}

\subsection{Experiment 2}
\begin{figure*}
\begin{tabular}{ccc}{\bf Vertical well} &{\bf Horizontal well}	&{\bf Inclined well} \\ 
\includegraphics[width=0.3\linewidth]{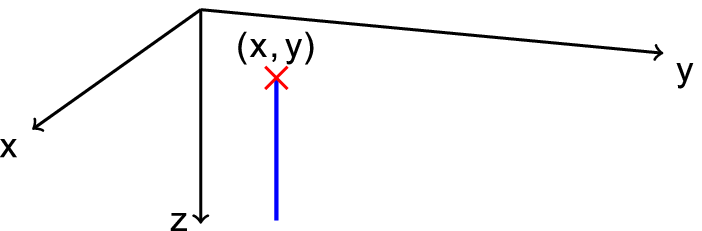} &\includegraphics[width=0.3\linewidth]{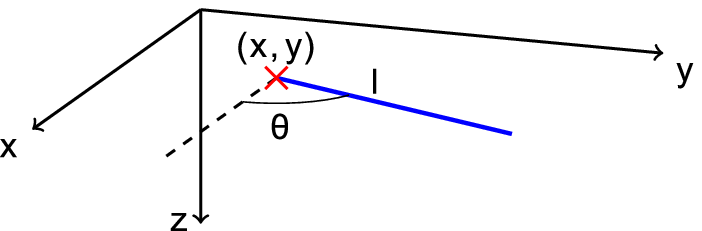} &\includegraphics[width=0.3\linewidth]{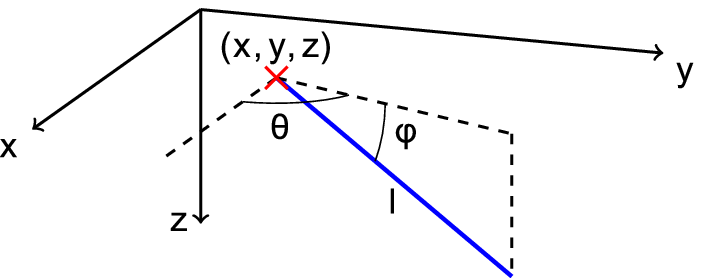}
\end{tabular}
\caption{Types of well and parameterizations used in Experiments 1, 2 and 3, respectively. Left: vertical well parameterized by $(x,y)$; centre: arbitrarily oriented horizontal well parameterized by $(x,y,l,\theta)$; right: inclined well parameterized by $(x,y,z,l,\theta,\phi)$. Main well bore is shown as a solid line, while the cross denotes the well heel.}\label{F:well_paths}
\end{figure*}

This experiment uses the same reservoir model, fluid properties, production period, BHP bounds and economic parameters as Experiment 1. As before, we consider a case with no production constraints (Case 2A) and one with maximum flowrate constraints (Case 2B); in this case, 1500 m$^3$/day for both injectors and producers. The key difference from Experiment~1 is that instead of drilling six vertical wells, we drill four arbitrarily oriented horizontal wells; two injectors and two producers. The positions of each of these wells are parameterized by four variables; $x$, $y$, $l$, and $\theta$. The co-ordinates of the well heel are given by $x$ and $y$, $l$ is the length of the well, and $\theta$ is the angle of the well in the $x$-$y$ plane, where an angle of zero corresponds to being oriented in the positive $x$ direction. This parameterization is illustrated in \fref{F:well_paths} (middle image).

The well length $l$ has a minimum value of 100 m and a maximum value of 320 m. The cost of drilling a longer well is incorporated into the NPV calculation by including a drilling cost of \$50,000 per metre of length, in addition to the base cost of \$25M per well. Thus a horizontal well may cost anywhere from \$30M to \$41M to drill. Given that the angle $\theta$ can take any value between 0 to 360$^{\circ}$, the horizontal wells will not, in general, be aligned with the grid, and so the well bores will not pass through the centres of the grid cells. As such, it makes sense for $x$ and $y$ to be specified in terms of the position of the well heel in metres, rather than as grid indices. One can then determine the length of intersection of the well bore with each cell through which it passes, which is used by IMEX to compute an appropriate well index for each segment.

\subsection{Experiment 3}
\begin{figure*}
\includegraphics[width=0.9\linewidth]{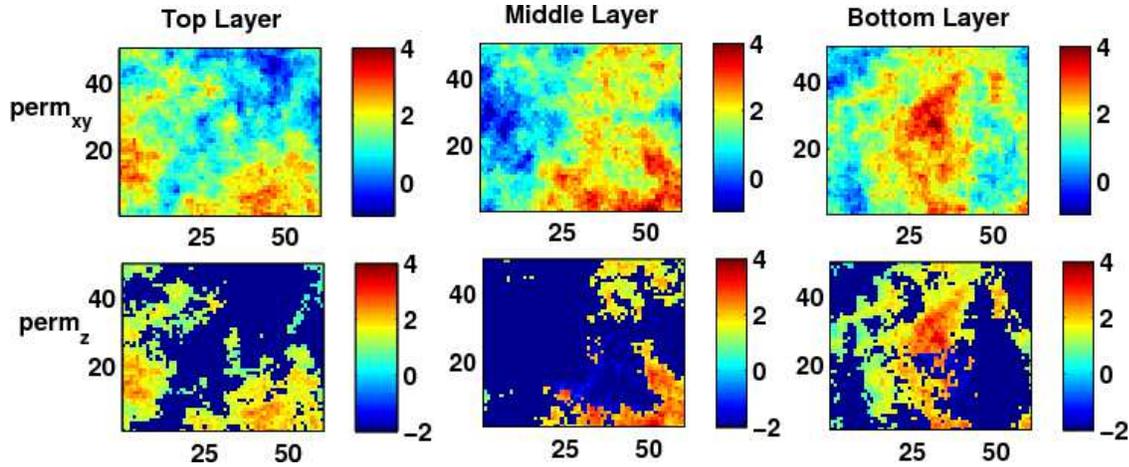}
\caption{Permeability field for the three-layer reservoir used in Experiment 3. Top row shows permeability (mD, log scale) in the $x$-$y$ directions, bottom row in the $z$ direction. Taken from the~\citet{SPE10}.}\label{F:perm3}
\end{figure*}

In this experiment, we consider the problem of drilling inclined wells in three dimensions. Each well's position can now be parameterized by six variables; the same $x$, $y$, $l$ and $\theta$ as in Experiment 2, as well as $z$ (the depth of the well heel) and $\phi$ (angle the well makes with the horizontal plane). We use a production period of 15 years, with well control periods of 3 years. Thus each well is parameterized by six positional parameters and five control parameters, for a total of eleven variables. The same bounds on BHP for injectors and producers are used as in the previous two experiments, as well as the same economic parameters. The synthetic reservoir is a 60$\times$50$\times$3-cell grid with cell dimensions 40$\times$40$\times$25 m (total field size: 2400$\times$2000$\times$75 m) and a depth of 2000 m. The permeability field of the reservoir (shown in \fref{F:perm3}) is isotropic in the $x$ and $y$ directions, with a much lower average permeability in the $z$ direction. The depth to water-oil contact is 2060 m, and so the top two layers of the reservoir initially contain mostly oil, and the bottom layer mostly water. 

We consider the problem of placing two injection and two production wells (44 variables total). The same NPV calculation is used as in Experiment~2, including the additional drilling cost per unit length. The length of each well must be at least 100 m and at most 400 m. To reduce the number of infeasible points in solution space, we place the following restrictions on the positional parameters:

\begin{enumerate}
\item The angle with the vertical, $\phi$, of each well must lie between 0 and 10 degrees. Given that the minimum well length is 100 m, steeper drilling angles will tend produce wells with toes that lie below the bottom layer of the reservoir, which correspond to invalid input.
\item The $z$ value of the well heel for injectors must lie between 25 to 50 m, while for producers it must lie between 0 to 50 m. This means that the heel of an injection well must lie in the middle layer of the reservoir, while that of a producer well can lie in the top or middle layer. As with the first restriction, this helps to limit the number of well configurations that are invalid due to the well toe position. It is also a sensible constraint since good solutions typically involve drilling producers in upper layers of the reservoir and injectors in lower layers. 
\end{enumerate}

As before, we consider two cases: one with no production constraints (Case 3A) and one with maximum fluid injection and production constraints of 5000 m$^3$/day at each well (Case 3B).

\section{Results} \label{S:results}
\begin{figure*}
\begin{tabular}{ccc}
{\bf Case 1A} &{\bf Case 2A} &{\bf Case 3A}\\
\includegraphics[width=0.32\linewidth]{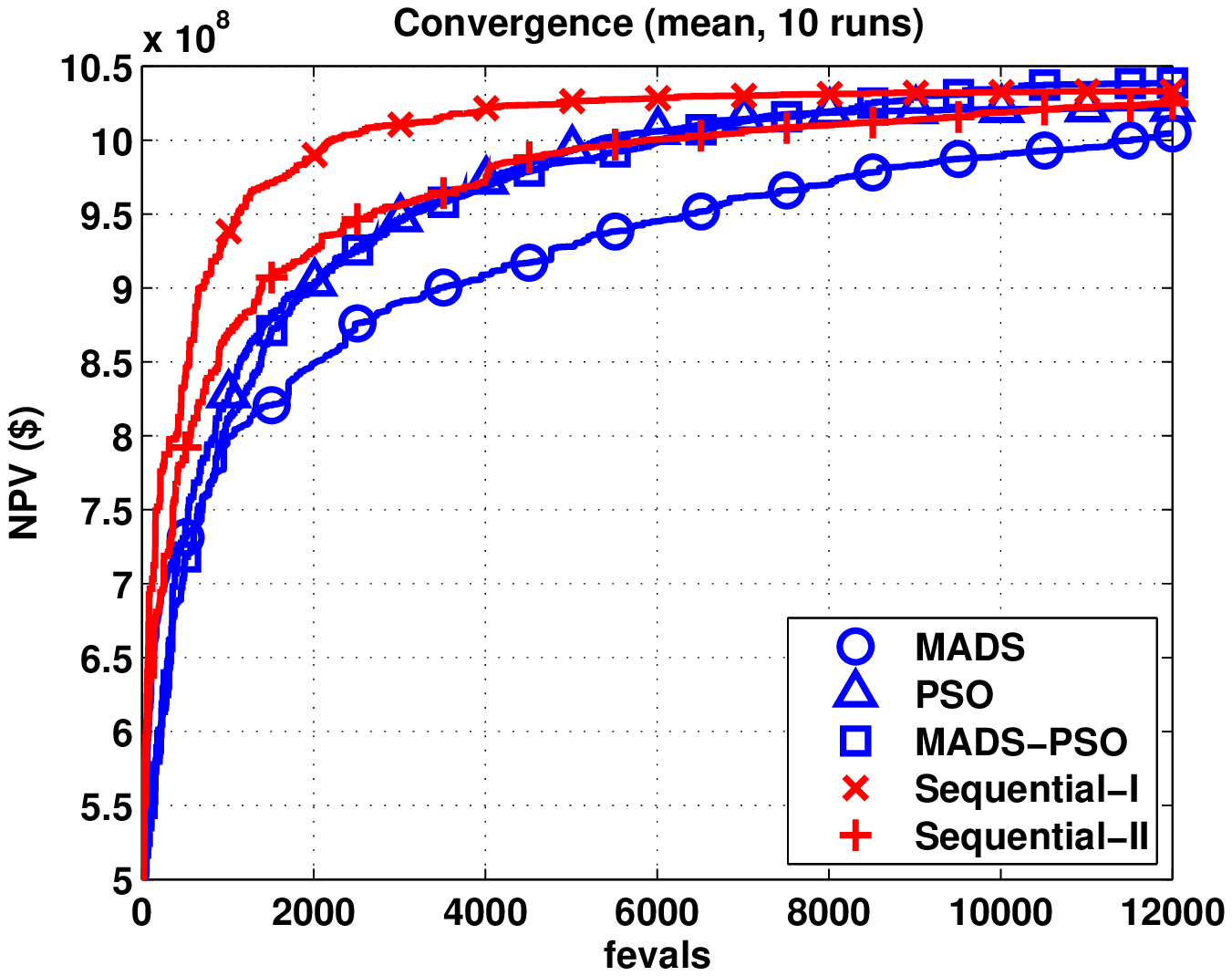} &\includegraphics[width=0.32\linewidth]{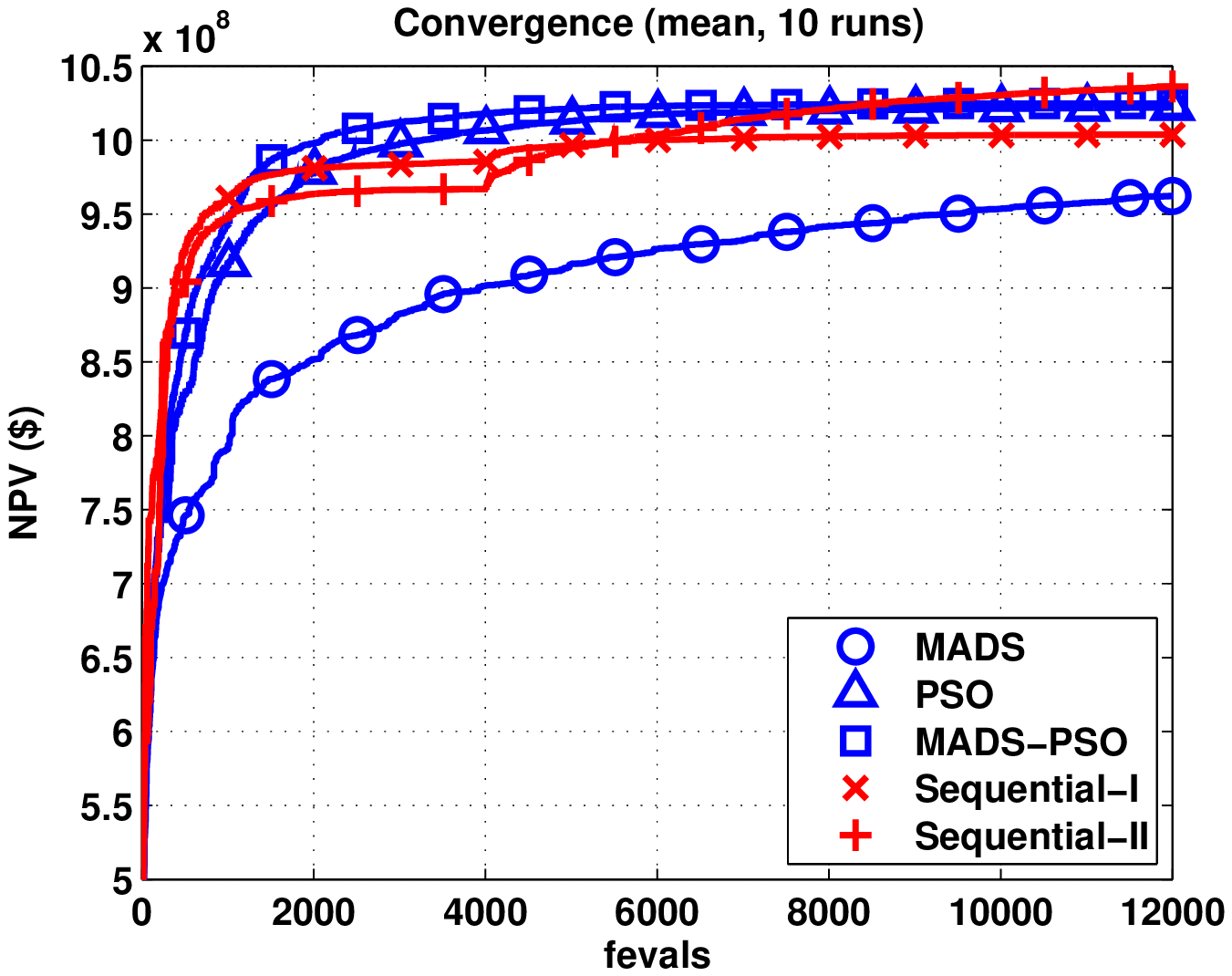} &\includegraphics[width=0.32\linewidth]{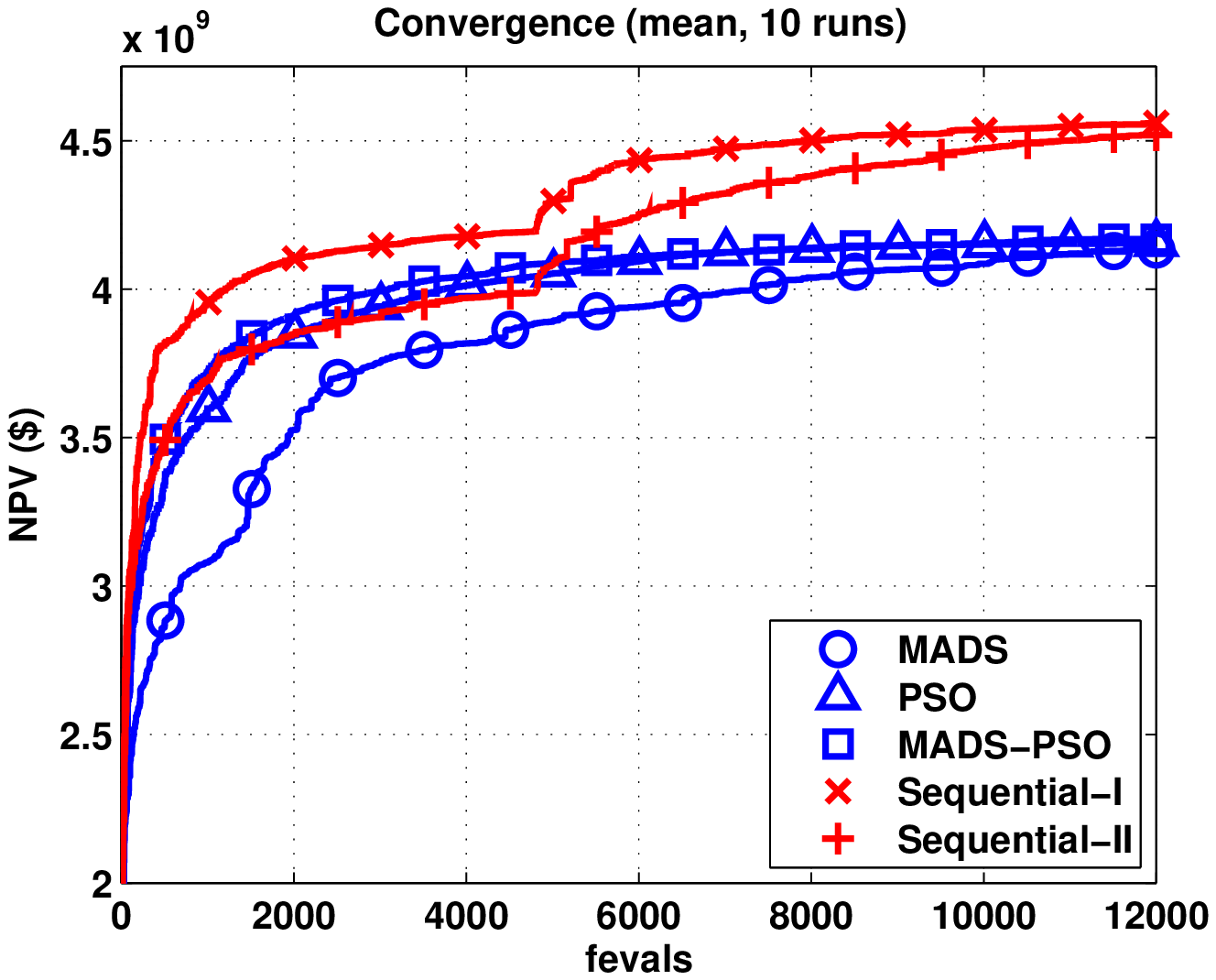} \smallskip \\
{\bf Case 1B} &{\bf Case 2B} &{\bf Case 3B}\\
\includegraphics[width=0.32\linewidth]{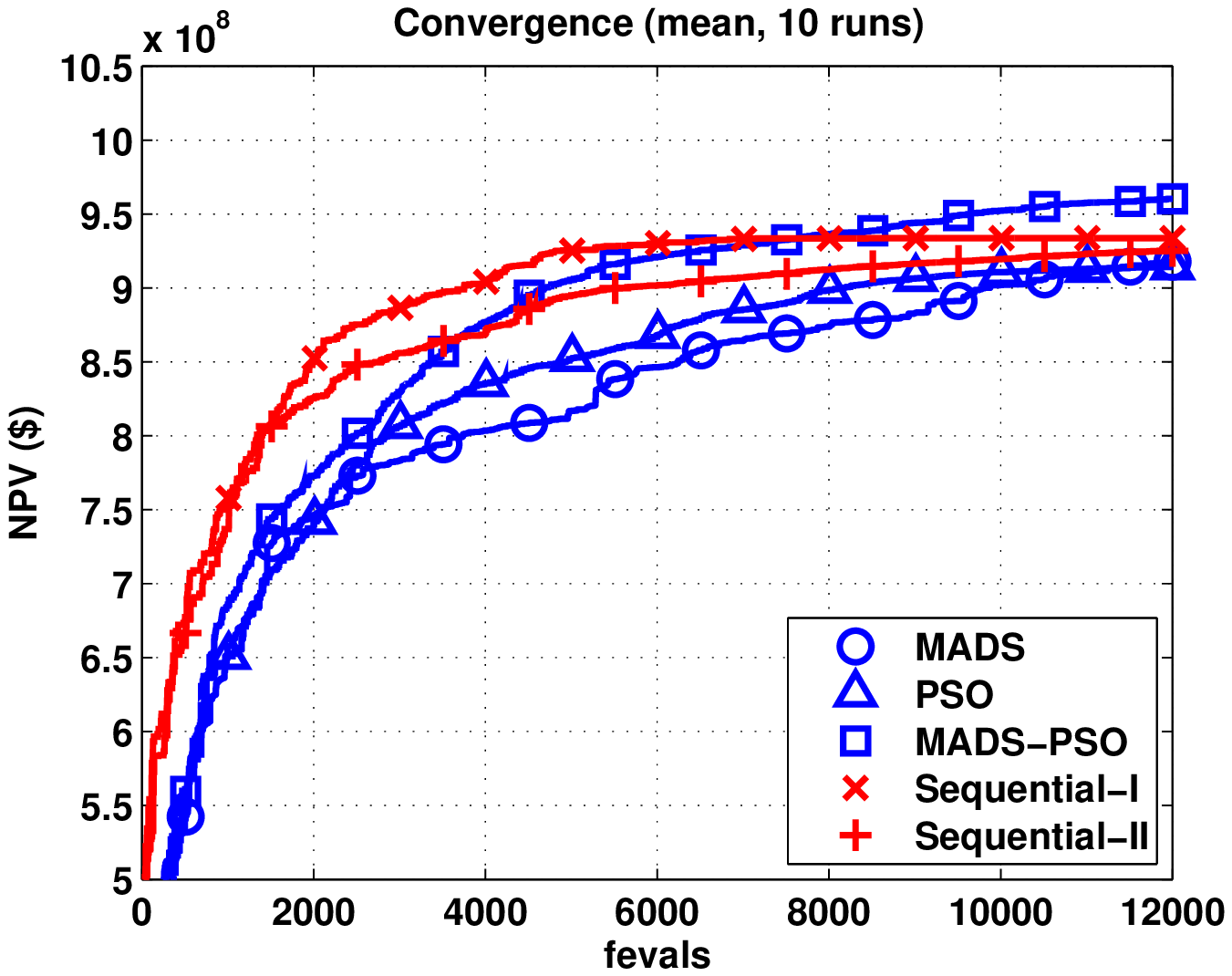} &\includegraphics[width=0.32\linewidth]{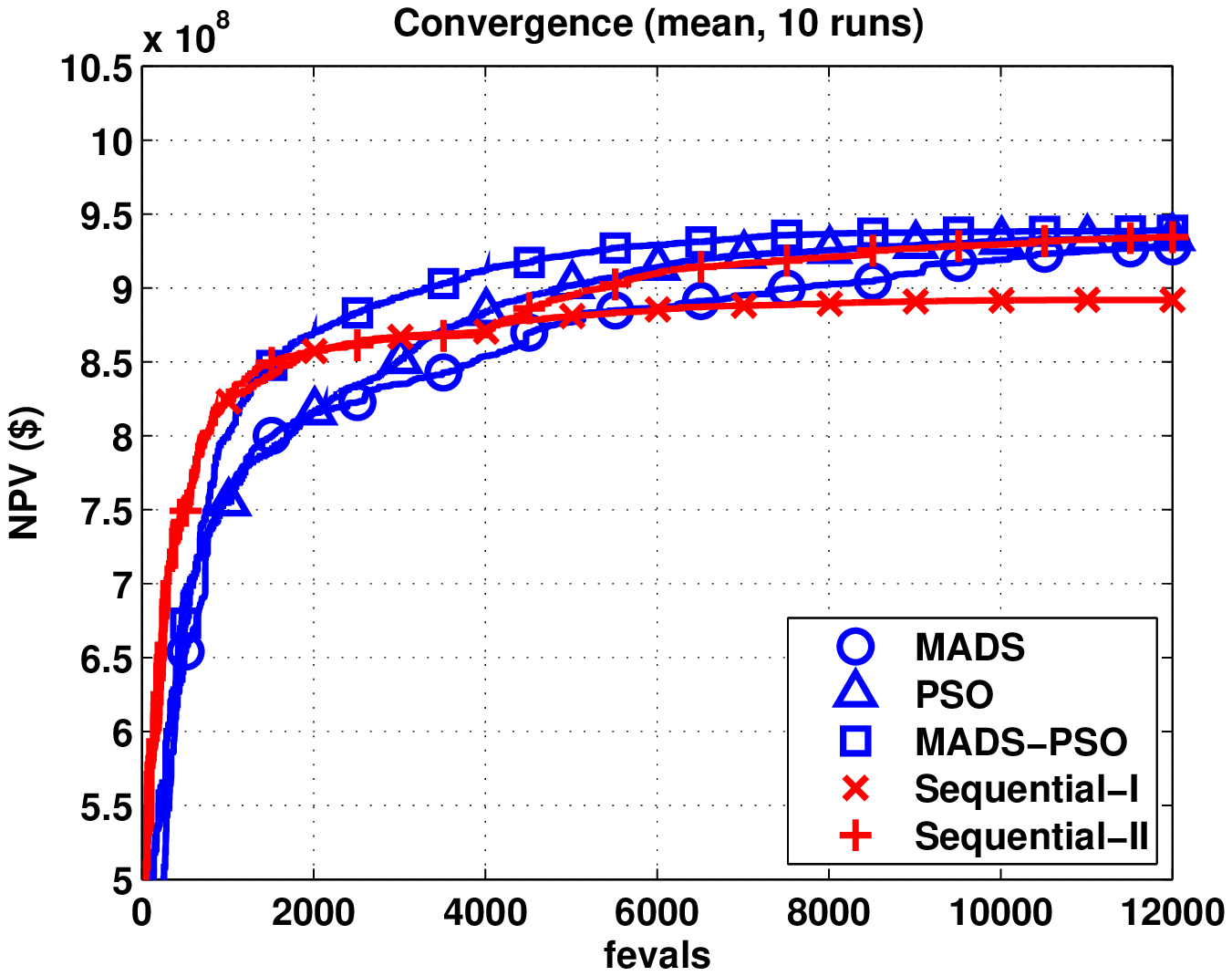} &\includegraphics[width=0.32\linewidth]{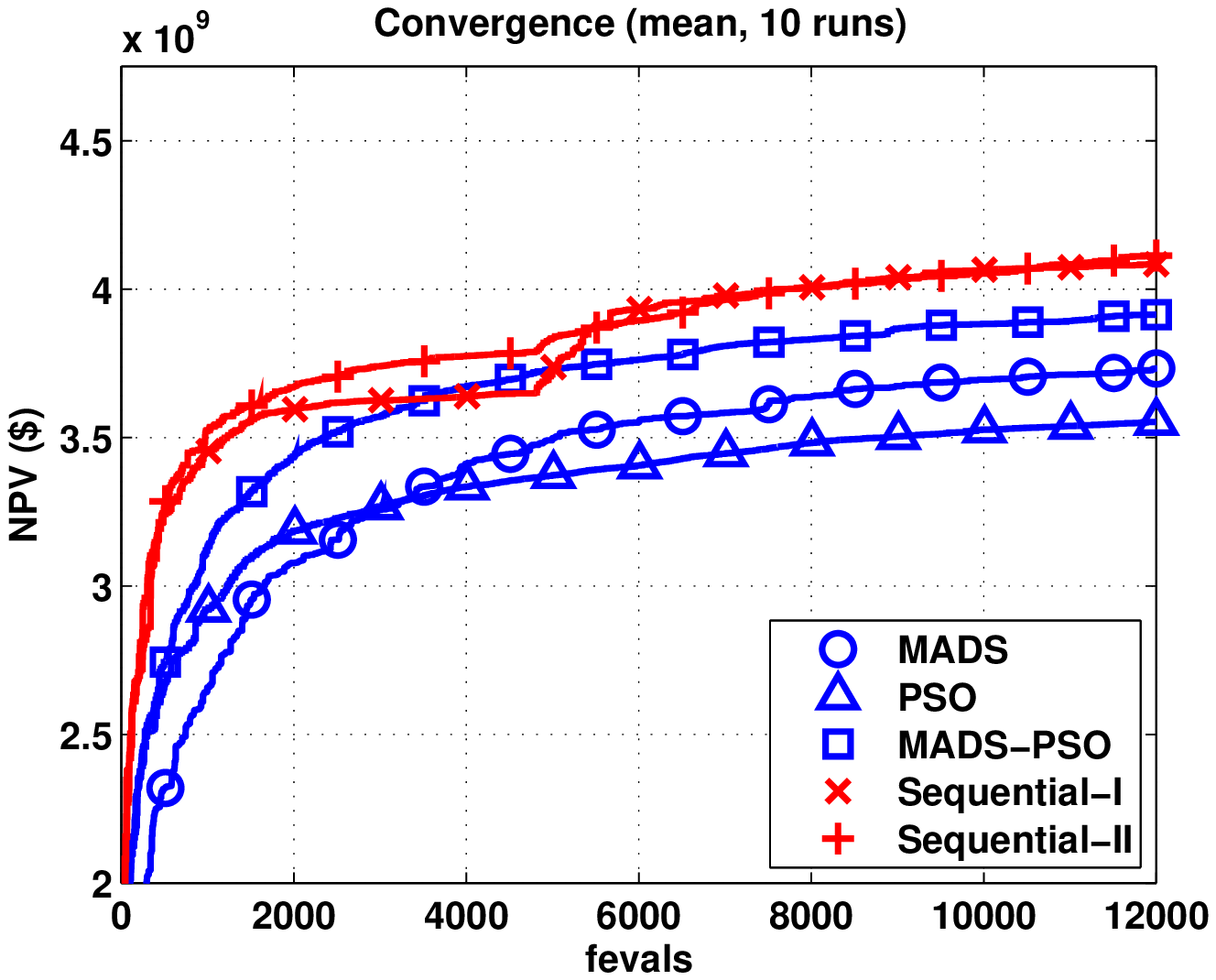}
\end{tabular}
\caption{Convergence plots for all experiments, showing NPV averaged over all 10 runs of each method, as a function of the number of simulations (fevals). Note that the y-axis scale for Cases 3A and 3B is different.}\label{F:convergence}
\end{figure*}

For each experiment, we tested five optimization approaches. The first three are simultaneous approaches (approaches which act on a vector of all parameters simultaneously), while the last two are sequential approaches, which consist of two optimization steps.
\begin{enumerate}
\item The MADS algorithm as described in Section~\ref{SS:MADS}, applied to all variables simultaneously. The algorithm was initialized by evaluating the objective function at 60 points using a Latin Hypercube search that is built in to NOMAD. The best point of these was selected as the initial incumbent point. The algorithm was run up to a maximum of 12,000 function evaluations, or until the mesh reached a specified minimum size.
\item The PSO algorithm as described in Section~\ref{SS:PSO}, implemented as a search class in NOMAD, with polling disabled. We ran the algorithm up to a maximum of 12,000 function evaluations, or until 100 consecutive iterations failed to improve the solution.
\item The MADS-PSO algorithm as described in Section~\ref{SS:MADS-PSO}. The algorithm was run up to a maximum of 12,000 function evaluations, or until the mesh reached a specified minimum size.
\item The Sequential-I approach as described in Section~\ref{SS:SEQ}. In Experiments 1 and 2, the first step consisted of running PSO for up to 4,000 function evaluations on the well placement optimization problem, and the second step then consisted of running MADS on the full problem for up to 8,000 function evaluations. In Experiment 3, we allocated 4,800 function evaluations to the first step and 7,200 to the second step, owing to the increased number of variables involved in Step 1. In all cases, injectors were held at their maximum BHP of 450 bar, and producers at the minimum BHP of 125 bar, during Step 1.
\item The Sequential-II approach as described in Section~\ref{SS:SEQ}. The procedure was the same as for Sequential-I, except that during Step 1, injectors were held at 425 bar, and producers at 150 bar.
\end{enumerate}
Every one of these approaches is nondeterministic, either as a result of the initialization procedure (for MADS) or the inclusion of a PSO component. Thus, to assess the overall performances of these algorithms, we ran each of them 10 times for every experiment. The results are summarized in Table~\ref{T:results}, which shows the best, worst, mean and standard deviation of the 10 runs for each algorithm. The best overall value for each experiment is highlighted in bold font. Convergence plots showing the NPV averaged over all 10 runs of each method as a function of the number of objective function evaluations (fevals) are shown in \fref{F:convergence}.

\begin{table}
\caption{Results for all experiments. Values shown are NPV in \$$\times 10^8$. Best values for each case are highlighted in bold.}\label{T:results}. 
\begin{tabular}{llrrrr}
{\bf Case}	&{\bf Algorithm}		&{\bf Best}		&{\bf Worst} 	&{\bf Mean}		&{\bf St. Dev} \\ 
 \hline \noalign{\smallskip}
1A		&{ MADS}			&10.53		&9.185		& 10.05		&0.548 \\
		&{ PSO}			&10.61		&9.316		& 10.22		&0.418 \\
		&{ MADS-PSO}		&10.97		&9.277		& {\bf 10.39}	& 0.534\\
		&{ Sequential-I}		&10.77		&{\bf 9.969}		& 10.33		&{\bf 0.306} \\
\smallskip	&{ Sequential-II}		&{\bf 10.99}		&9.003		& 10.25		&0.556 \\

1B 		&{ MADS}			&9.945		&8.241		&9.179		&{\bf 0.517} \\
		&{ PSO}			&9.927		&8.109		&9.143		&0.552 \\
		&{ MADS-PSO}		&{\bf 10.36}		&{\bf 8.768}		&{\bf 9.603}		&0.528 \\
		&{ Sequential-I}		&10.24		&8.160		&9.336 		&0.574 \\
\smallskip	&{ Sequential-II}		&10.12		&7.413		&9.251 		&0.844 \\

2A		&{ MADS}			&10.51		&7.974		&9.623		&0.783 \\
		&{ PSO}			&10.70		&9.406		&10.22		&0.538 \\
		&{ MADS-PSO}		&10.80		&9.652		&10.25		&0.439 \\
		&{ Sequential-I}		&10.47		&9.620		&10.04 		&0.316 \\
\smallskip	&{ Sequential-II}		&{\bf 10.83}		&{\bf 9.981}		&{\bf 10.36}		&{\bf 0.287} \\

2B		&{ MADS}			&10.20		&{\bf 8.861}		&9.282		&{\bf 0.441} \\
		&{ PSO}			&9.963		&8.198		&9.342		&0.475 \\
		&{ MADS-PSO}		&10.16		&8.217		&{\bf 9.390}		&0.597 \\
		&{ Sequential-I}		&9.864		&7.505		&8.920 		&0.701 \\
\smallskip	&{ Sequential-II}		&{\bf 10.32}		&8.410		&9.345		&0.668 \\

3A		&{ MADS}			&46.03		&34.03		&41.34		&4.69 \\
		&{ PSO}			&46.20		&32.69		&41.55		&4.47 \\
		&{ MADS-PSO}		&45.42		&33.83		&41.69		&3.62 \\
		&{ Sequential-I}		&{\bf 48.50}		&42.13		&{\bf 45.60} 	&2.40 \\
\smallskip	&{ Sequential-II}		&47.42		&{\bf 43.94}		&45.19		&{\bf 1.00} \\

3B		&{ MADS}			&45.50		&30.70		&37.33		&5.43 \\
		&{ PSO}			&41.83		&26.97		&35.53		&5.32 \\
		&{ MADS-PSO}		&45.19		&30.64		&39.15		&4.82 \\
		&{ Sequential-I}		&{\bf 46.14}		&31.57		&40.84 		&5.31 \\
		&{ Sequential-II}		&44.46		&{\bf 38.16}		&{\bf 41.13}		&{\bf 2.04} \\
\end{tabular}
\end{table}

\section{Discussion} \label{S:discussion}
We first comment on some trends that are apparent across all (or nearly all) of the experiments. With respect to the three simultaneous approaches, one clearly sees a benefit in hybridizing PSO and MADS for this problem, as compared to applying only MADS or only PSO. The average performance of MADS-PSO was better in all six experiments, although it did not always find the best solution overall between the three. The convergence plots in Fig.~\ref{F:convergence} indicate that MADS-PSO generally found good solutions more quickly than the other two simultaneous approaches as well. 

Additionally, we note that for the unconstrained problems (Cases 1A, 2A and 3A), PSO generally performed better than MADS on average, while in the constrained problems MADS performed just as well in Cases 1B and 2B, and significantly better in Case 3B. This may indicate that the nonlinear constraint handling strategy for PSO described in Section~\ref{SS:constr} is less effective than the strategy used by MADS, and could be improved upon. Ranking any feasible particle ahead of any infeasible particle drives the search towards feasible regions quickly and thus may disregard some promising infeasible solutions in other regions of the search space. Thus it may be worth considering more sophisticated constraint-handling methods that have been proposed for PSO \citep[e.g.][]{HW07b, KLZZ09}.

\begin{figure}
\includegraphics[width=0.75 \linewidth]{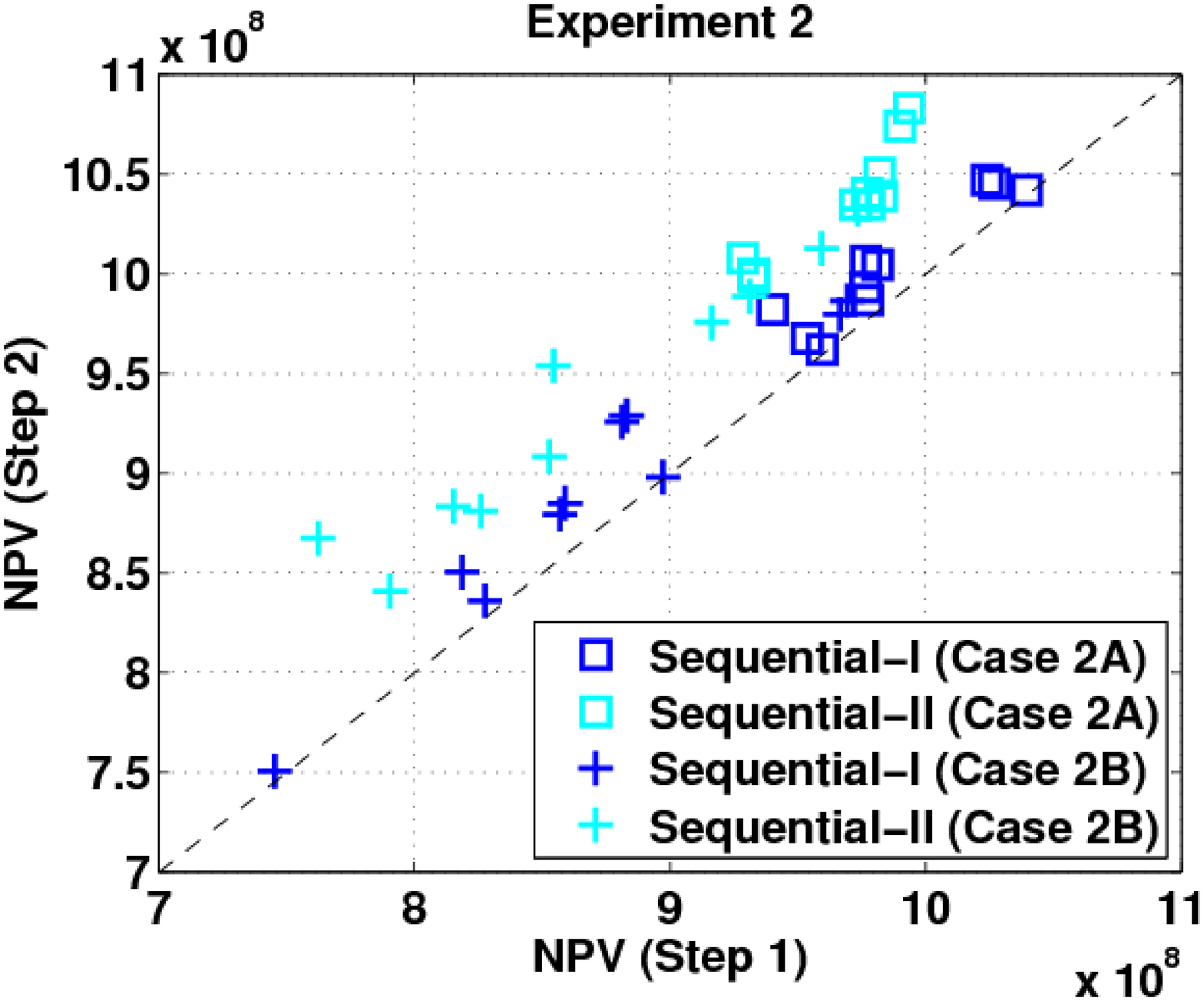}
\caption{Scatter plot showing results of Sequential-I and Sequential-II approaches for Cases 2A and 2B. Dashed line is the line $y=x$.}\label{F:exp2_scatter}
\end{figure}

The performance of Sequential-I versus Sequential-II appears to be somewhat context-specific. In Experiments 1 and 3, the average performance of both approaches was comparable (less than 1\% difference in mean NPV in all four cases). In Experiment 2, however, the Sequential-II technique provided results that were 4-5\% better than those of Sequential-I, on average, for both test cases. In \fref{F:exp2_scatter}, we show a scatter plot of the results of these two approaches for Cases 2A and 2B, in which the NPV found after Step 1 of the approach is plotted against the (final) NPV found after Step 2, for all 10 runs of each approach. All the points lie above the line $y=x$ since the second step cannot produce a solution with a lower NPV. It is clear from the plot that there were many instances where the best solution found by Step 1 of Sequential-I was not significantly improved by Step 2 (i.e., the points lying close to the line $y=x$). This suggests that the first step of Sequential-I tended to find a well configuration that was nearly optimal for the assumed control scheme (maximum/minimum BHP values at injectors and producers, respectively), and that it was difficult to improve on this solution in Step 2. The first step of Sequential-II, on the other hand, found solutions that were significantly improved during the second step, which eventually produced better results overall in this experiment. These results suggest that the Sequential-II approach is more robust than Sequential-I.

Of particular interest is the performance of the simultaneous MADS-PSO approach versus the two sequential approaches. In Case 1A, the performance of MADS-PSO was fairly comparable to that of {Sequential-I} on average (within 1\%), but in Case 1B, MADS-PSO outperformed both sequential approaches by roughly 3\% on average. In Experiment 2, Sequential-II and MADS-PSO had comparable average performances for both test cases. In Experiment 3, both sequential approaches provided much better results than MADS-PSO. The NPV of solutions found by Sequential-II, for instance, was 9\% higher on average in Case 3A, and 7\% higher in Case 3B. Thus, there appears to be a trend -- as the number of variables parameterizing each well's position increased (from two in Experiment~1 to six in Experiment~3), the performance of the sequential approaches relative to the simultaneous approaches improved.

\begin{figure}
\begin{tabular}{C{0.45\linewidth}C{0.45\linewidth}}
	\includegraphics[width=0.8\linewidth]{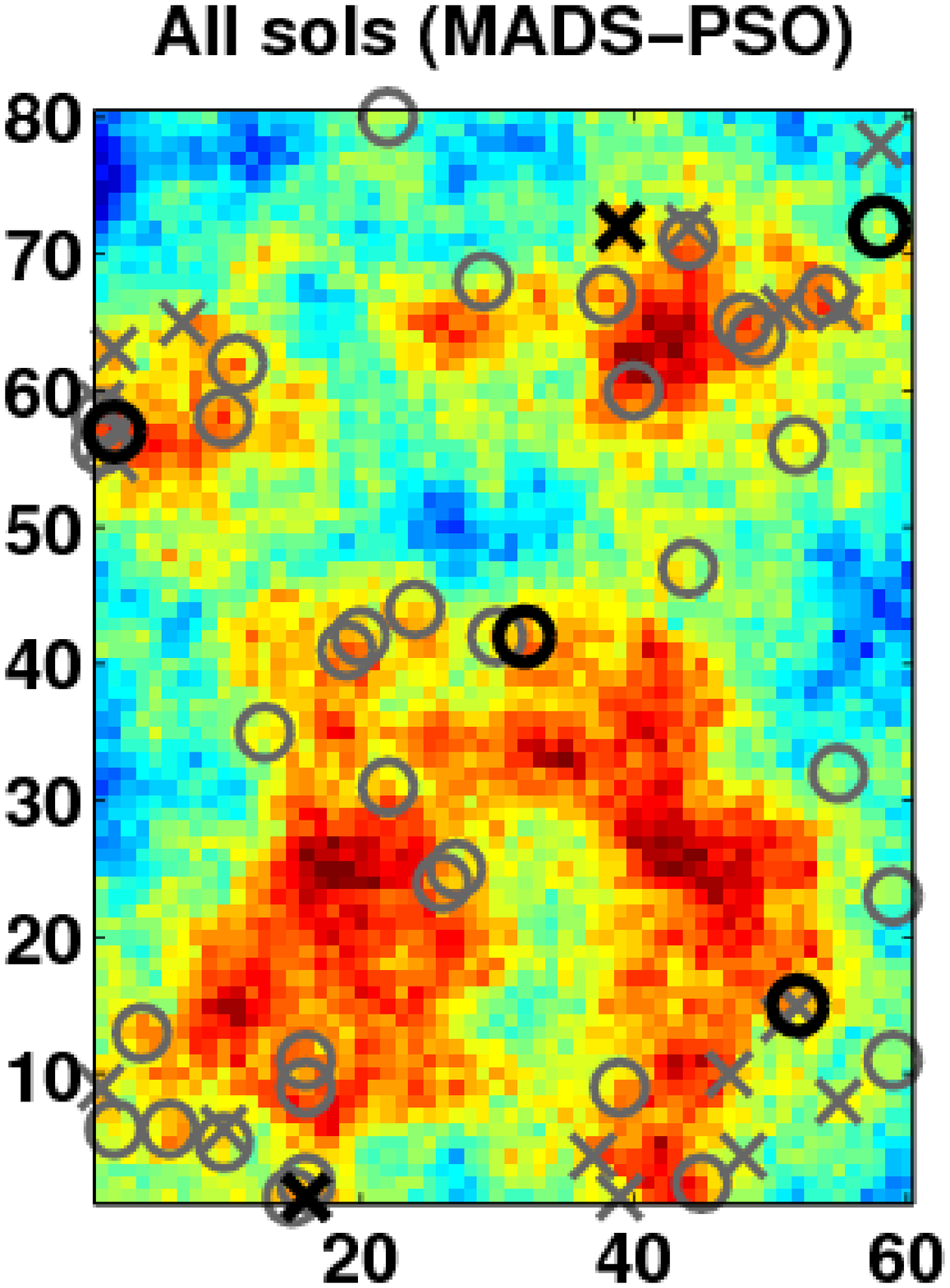}
&	\includegraphics[width=0.8\linewidth]{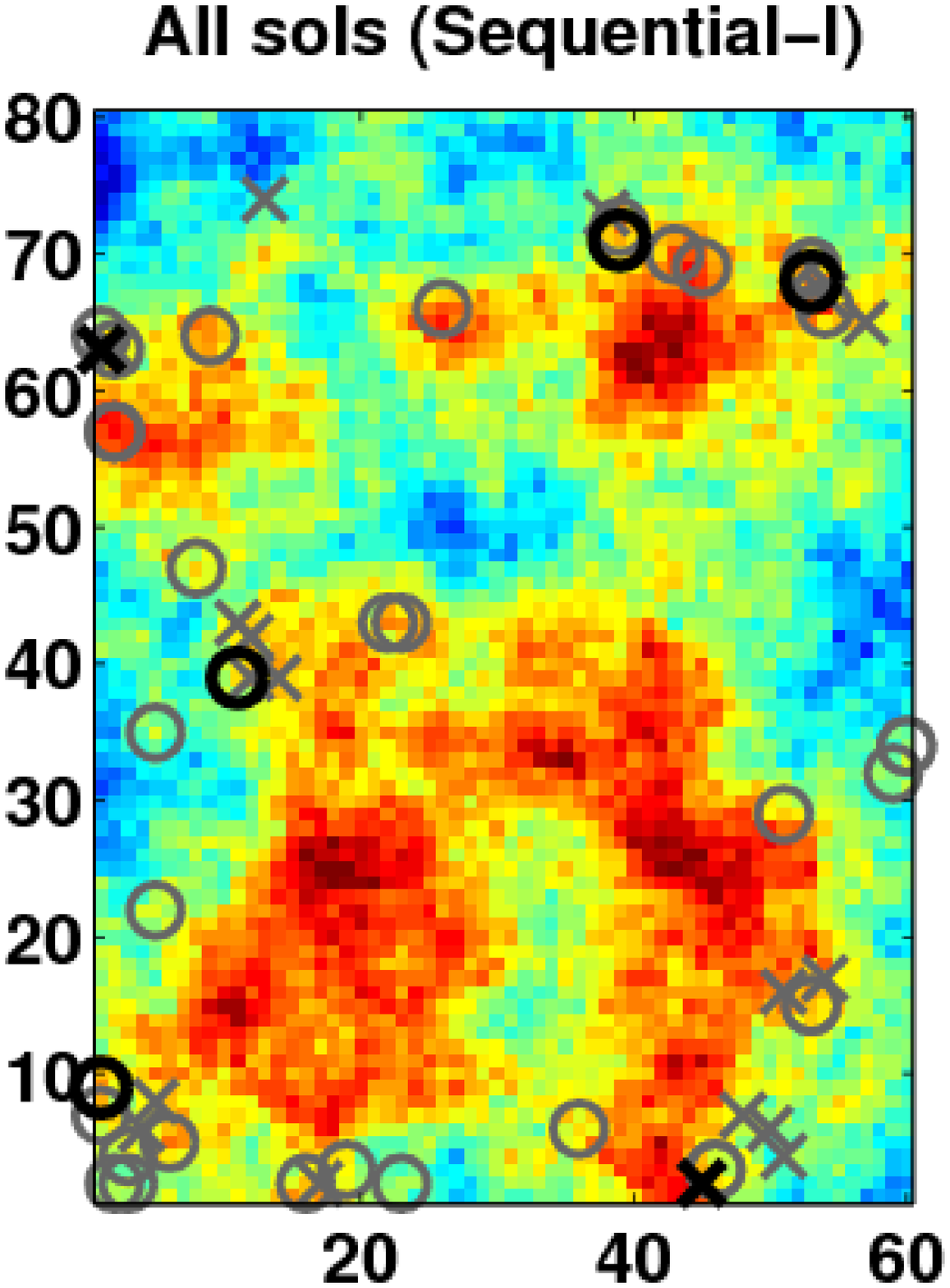}
\end{tabular}
\caption{Well positions corresponding to optimal solutions found at the end of each of the 10 runs of MADS-PSO (left) and Sequential-I (right), for Case 1B. Placement of injectors ($\times$) and producers ($\circ$) are overlaid on the permeability field (log scale, cf. \fref{F:perm}). The positions corresponding to the best overall solution found by each approach are highlighted with thick black symbols.}\label{F:Case1B_allsols}
\end{figure}

\begin{figure}
\begin{tabular}{C{0.45\linewidth}C{0.45\linewidth}}
	\includegraphics[width=0.9\linewidth]{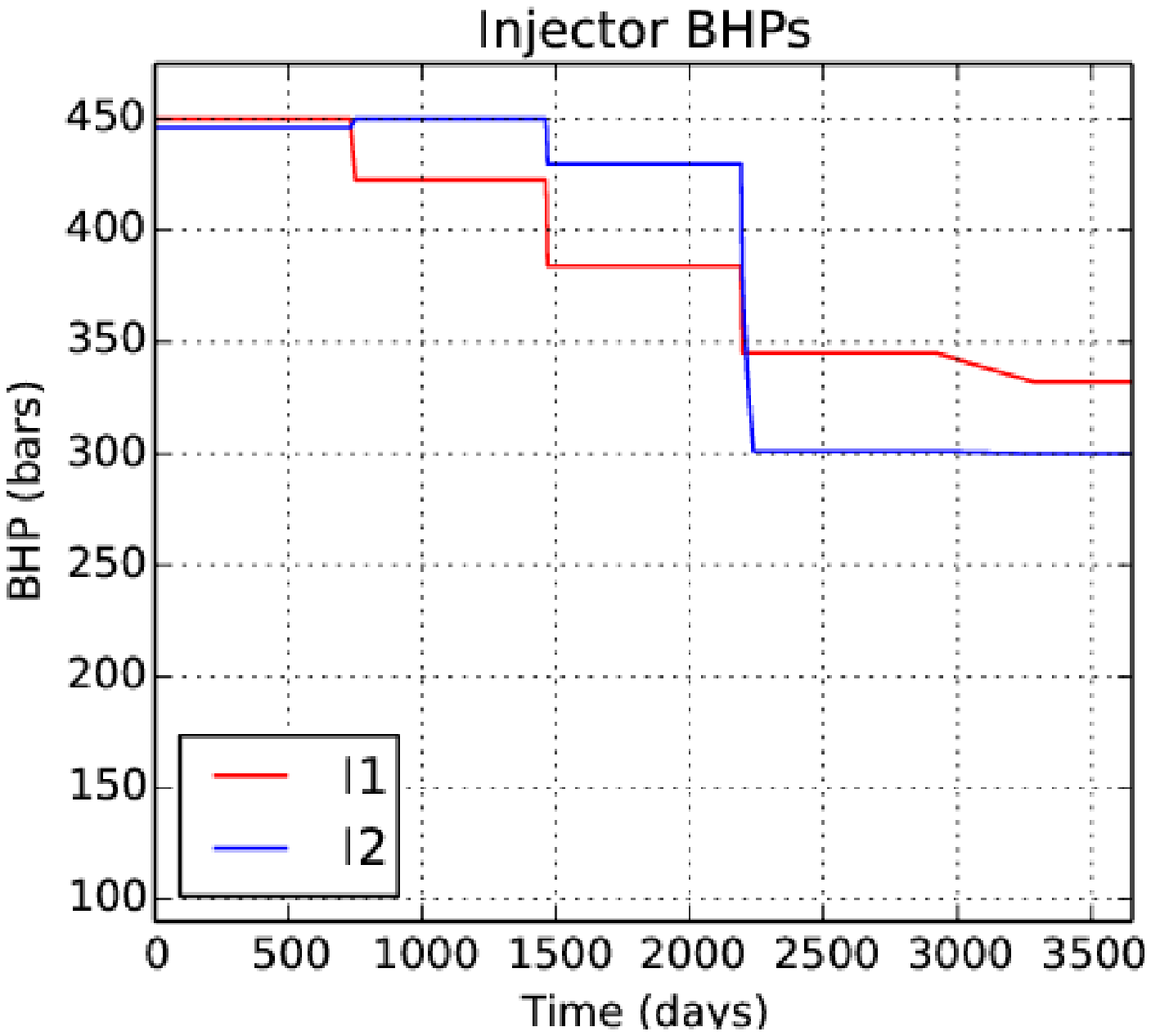}	 &\includegraphics[width=0.9\linewidth]{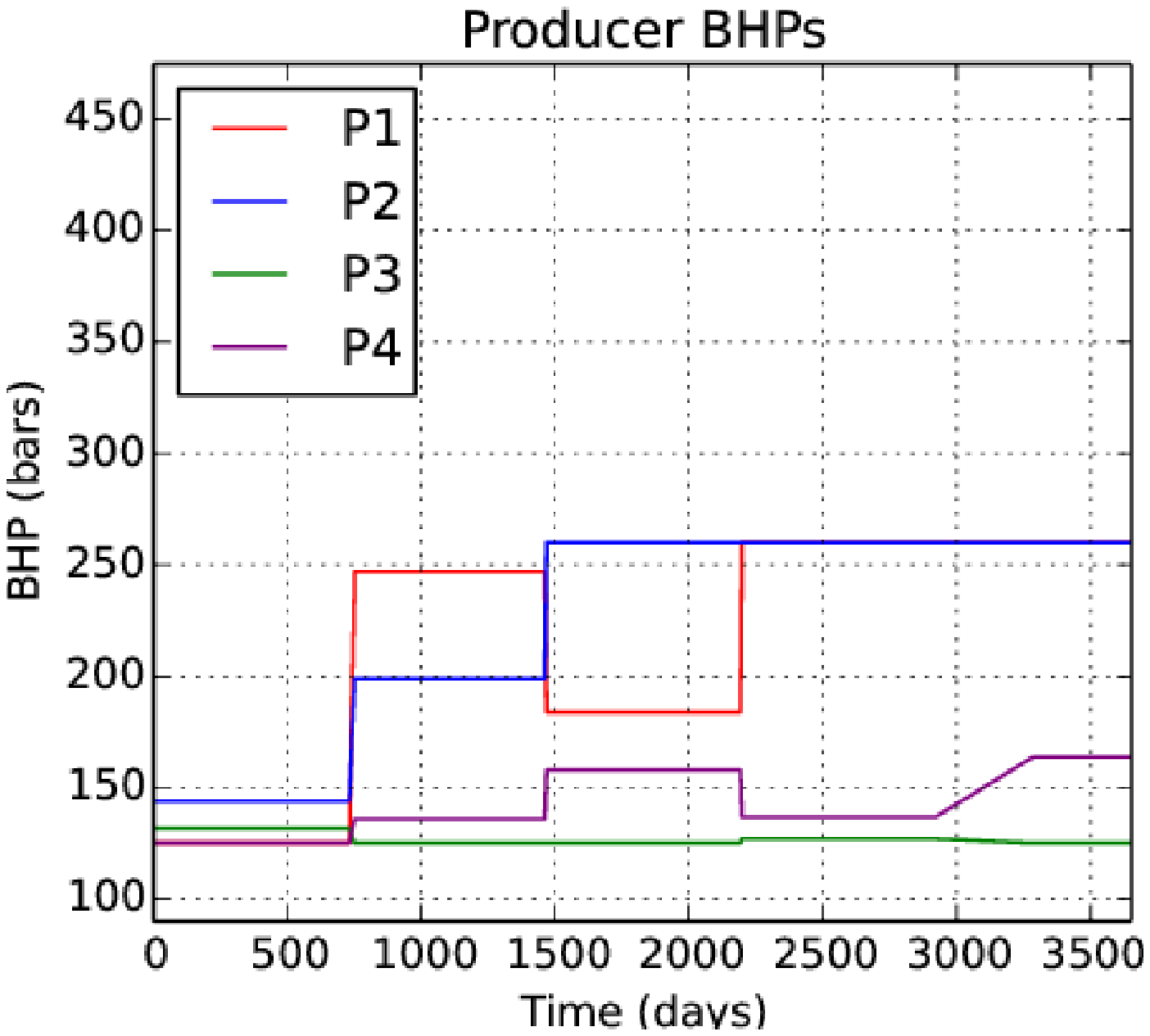} \\
	\includegraphics[width=\linewidth]{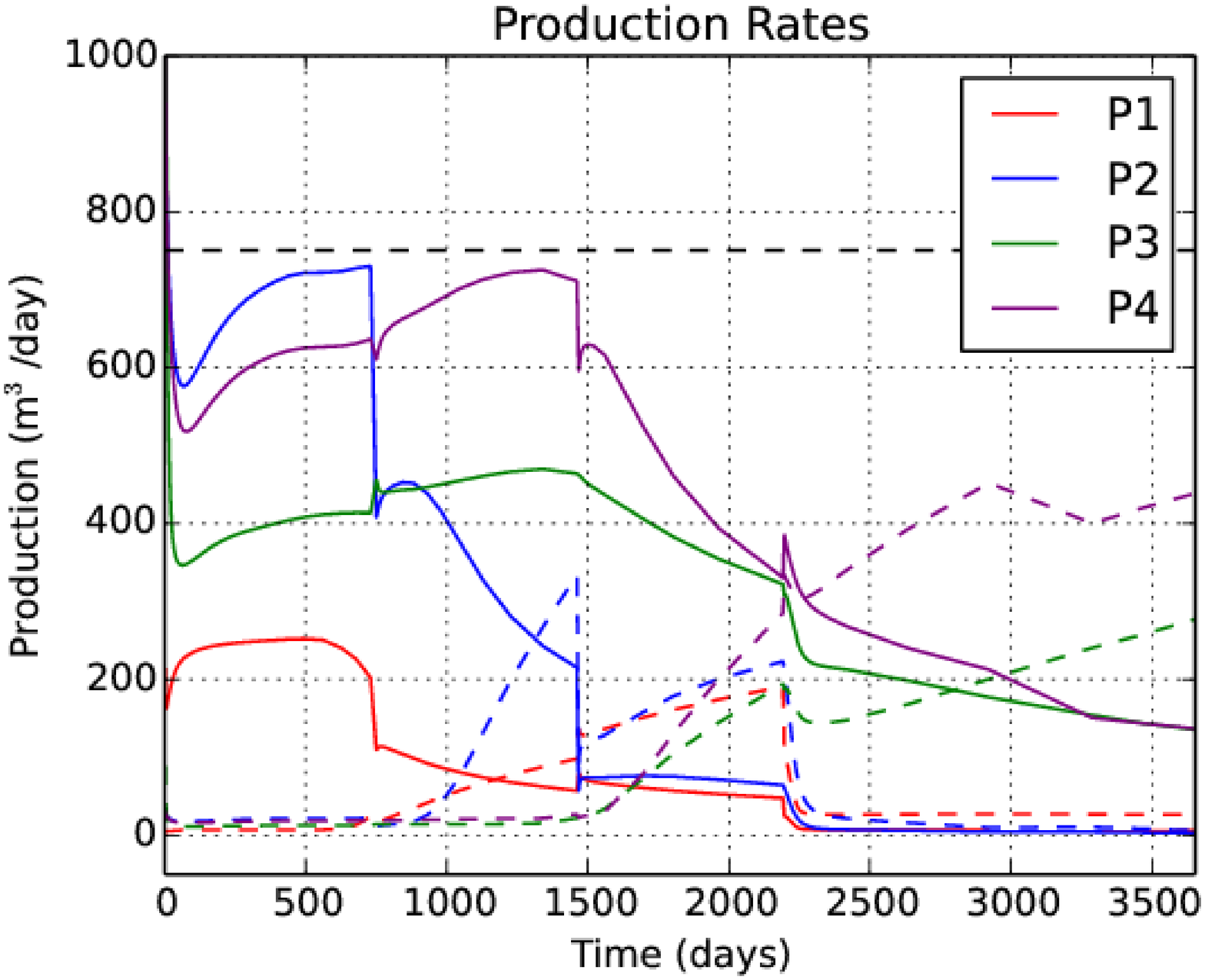} &	\includegraphics[width=\linewidth]{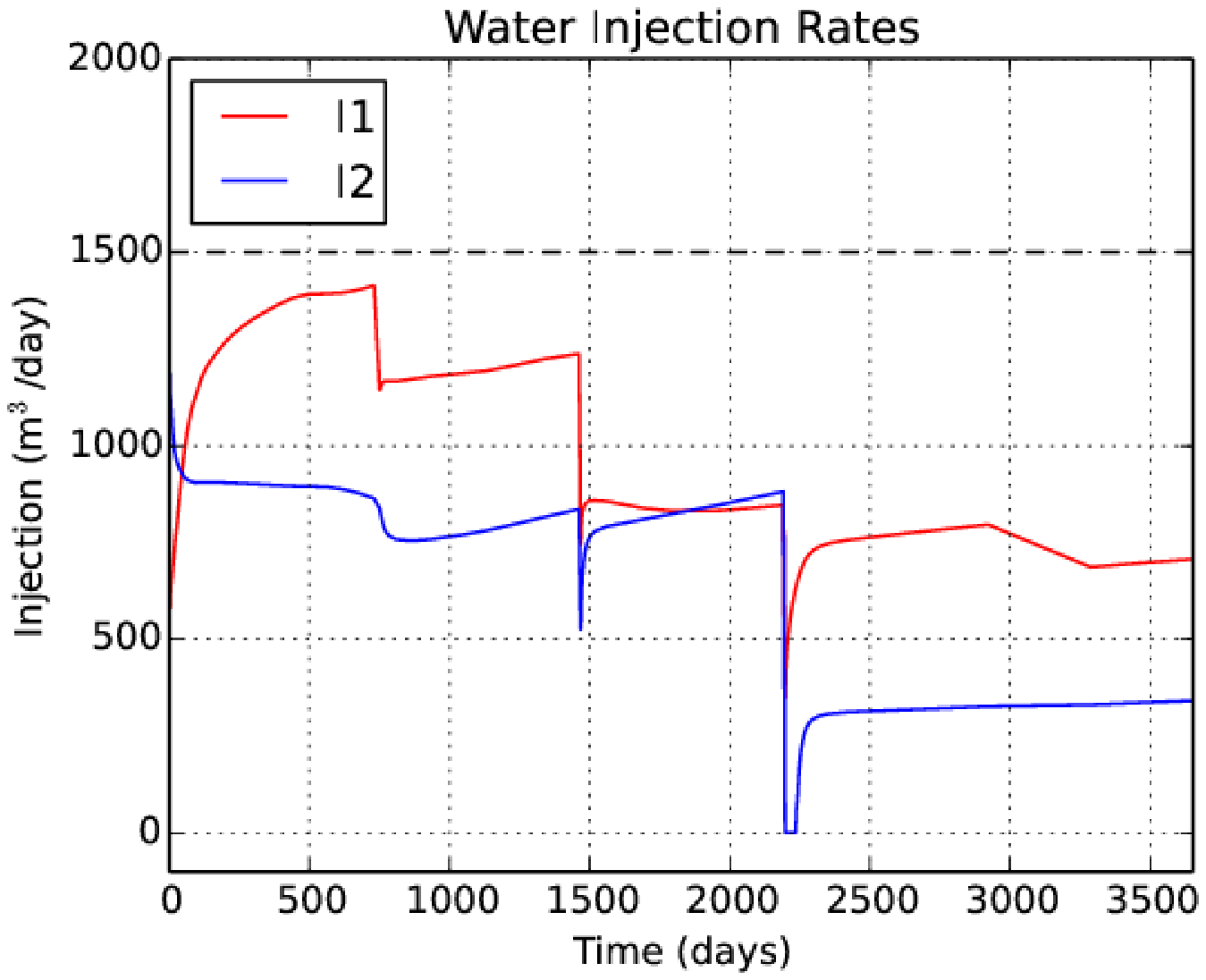} \\
\end{tabular}
\caption{BHP controls and corresponding flow rates for the best solution found by MADS-PSO for Case 1B, with an NPV of \$1.036$\times 10^{9}$. BHP controls for each injector and producer are shown in top row. Fluid production rates for each producer are shown in bottom left plot and water injection rates for each injector in bottom right. The solid lines in the bottom left plot indicate oil production; dashed lines indicate water production. The horizontal dashed lines indicate the maximum rate constraints for production and injection.}\label{F:Case1B_ctrls}
\end{figure}

We focus on the two cases where we observed the largest disparities in performance between simultaneous MADS-PSO and the sequential approaches. In Case 1B, where MADS-PSO gave better results, the goal was to place six vertical wells under the production constraints described in Section~\ref{SS:exp1}. In \fref{F:Case1B_allsols} we show the well positions corresponding to the best solution found by each of the 10 runs of the MADS-PSO (left figure) and Sequential-I (right figure) algorithms. We note that the solutions found by Sequential-I tend to place more wells towards the edges of the field and fewer towards the centre, compared to the solutions found by MADS-PSO. This occurs as a consequence of the control scheme used in the first step of Sequential-I. Holding injectors and producers at their maximum and minimum BHP values (respectively) generates high flow rates if wells are placed close together or in regions of high permeability, which result in constraint violations and lowered NPV due to premature waterflooding. Indeed, the overall best configuration of wells found by MADS-PSO (shown as thick black symbols in the left image of \fref{F:Case1B_allsols}) produces a constraint violation under the control schemes assumed in Step 1 of both the Sequential-I and Sequential-II approaches. The control scheme corresponding to this solution is shown in \fref{F:Case1B_ctrls}. It is clear that satisfying the constraints on maximum flowrate (1500 m$^3$/day for injectors and 750 m$^3$/day for producers) requires raising and/or lowering well BHPs several times (e.g. for the injector labeled I1 and producer labeled P2 after two years). Thus the first step of the two sequential approaches is not able to find some of the promising well configurations identified by the simultaneous approach. Although the second step of the sequential approaches allows the well positions to be altered in addition to the controls, it is primarily a local search and therefore is not likely to find a well configuration that is significantly different from the starting point provided by Step 1. 

%  A possible explanation is that satisfying the maximum flowrate constraints during Step 1 of the Sequential-I approach required placing wells further apart, since the well BHPs were held at their maximum/minimum values, which generates higher flow rates. The simultaneous MADS-PSO approach, on the other hand, could place wells in these regions since the BHPs could be raised or lowered to avoid violating a constraint. Thus the MADS-PSO approach was able to determine some well configurations which could not be found by the first step of the Sequential-I approach.  {\bf NB-- check well positions for best MADS-PSO solution with max/min controls and see if it produces a constraint violation}.

We note that \fref{F:Case1B_allsols} indicates that there is a great deal of variance in terms of the best positions found after each run. This is typical of problems involving well placement, since the objective function tends to contain many local optima with similar NPVs. Although this makes it harder to find the overall global optimum, it can be advantageous to have multiple different solutions with nearly optimal values, since some ``optimal'' solutions may be impractical for logistical reasons that are difficult to incorporate into the optimization algorithm (e.g. unsafe drilling conditions).

In Case 3A, where the two sequential approaches gave better results than MADS-PSO, the problem involved placing four inclined wells without any production constraints. The convergence plot (\fref{F:convergence}, top right) shows that even the solutions found during Step 1 of the Sequential-I approach, when well controls were held fixed, were better on average than those found by MADS-PSO. Thus, by focusing only on the 24 positional variables (rather than the 44 variables involved in the full placement-control problem), the first step of Sequential-I was quickly able to determine good well configurations. The solutions were then significantly improved during the second step, as is evident from the distinct bump in the convergence plots that occurs after 4800 fevals. In principle, the simultaneous MADS-PSO algorithm should be able to find these solutions as well, since the search space explored by the first step of Sequential-I is contained within the search space of the full simultaneous problem. However the convergence results would seem to indicate that the inclusion of the control variables along with the positional variables makes it more difficult for the simultaneous algorithm to determine optimal positional parameters, which play the most crucial role in the quality of the solution. 

\begin{figure*}
\begin{tabular}{C{0.32\linewidth}C{0.32\linewidth}C{0.32\linewidth}}
{\bf Top Layer} &{\bf Middle Layer} &{\bf Bottom Layer} \\
\includegraphics[width=\linewidth]{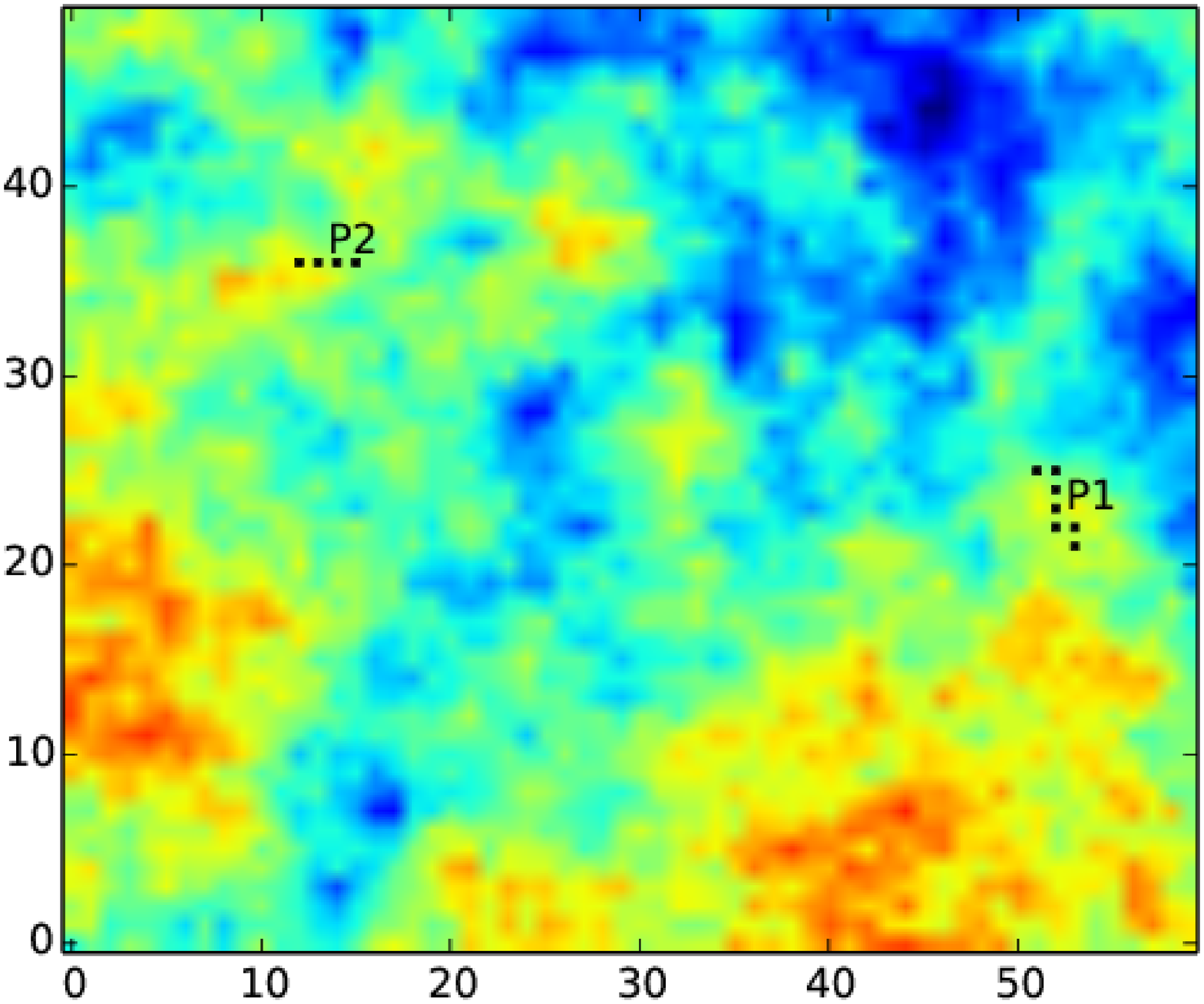}
&\includegraphics[width=\linewidth]{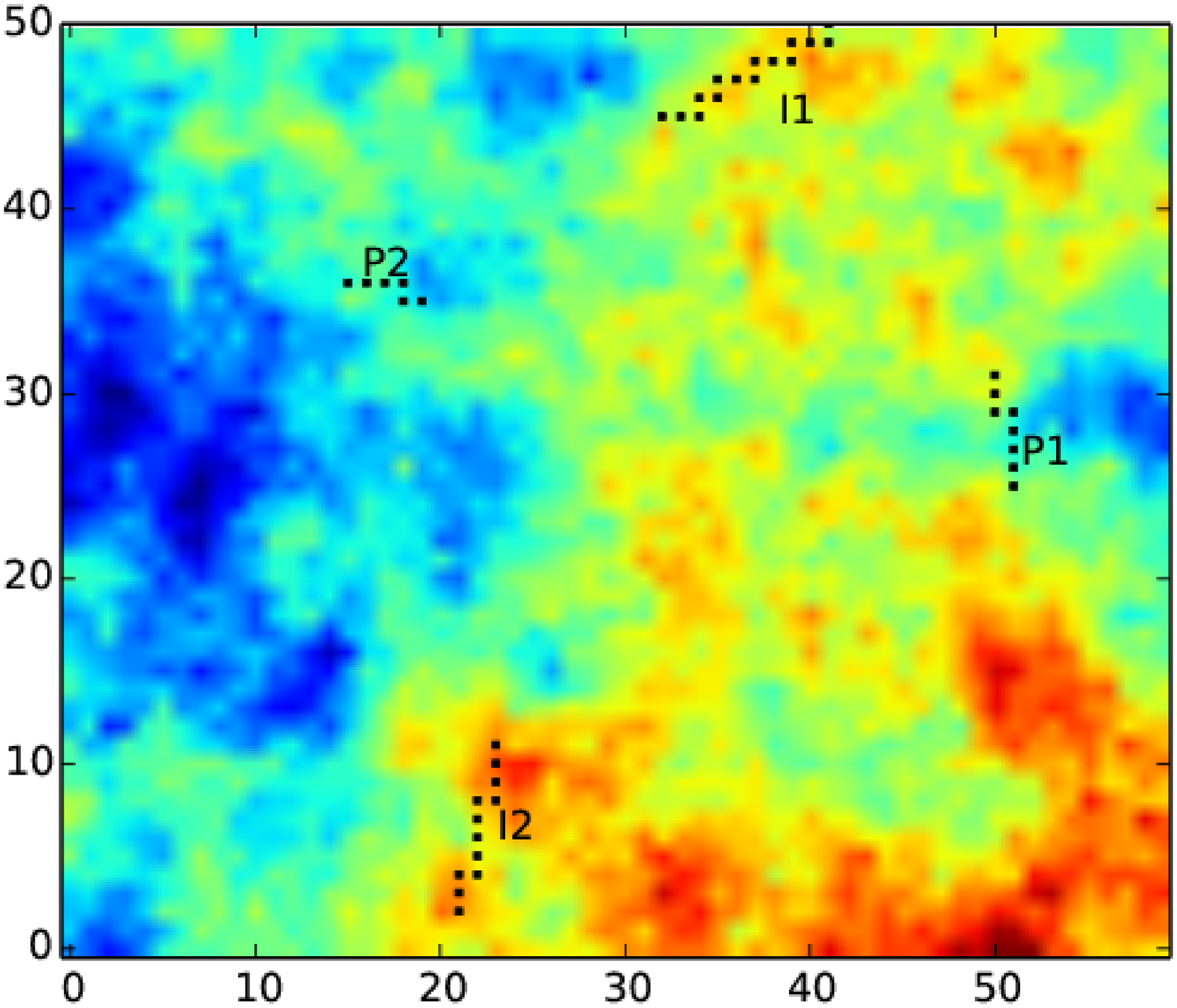}
&\includegraphics[width=\linewidth]{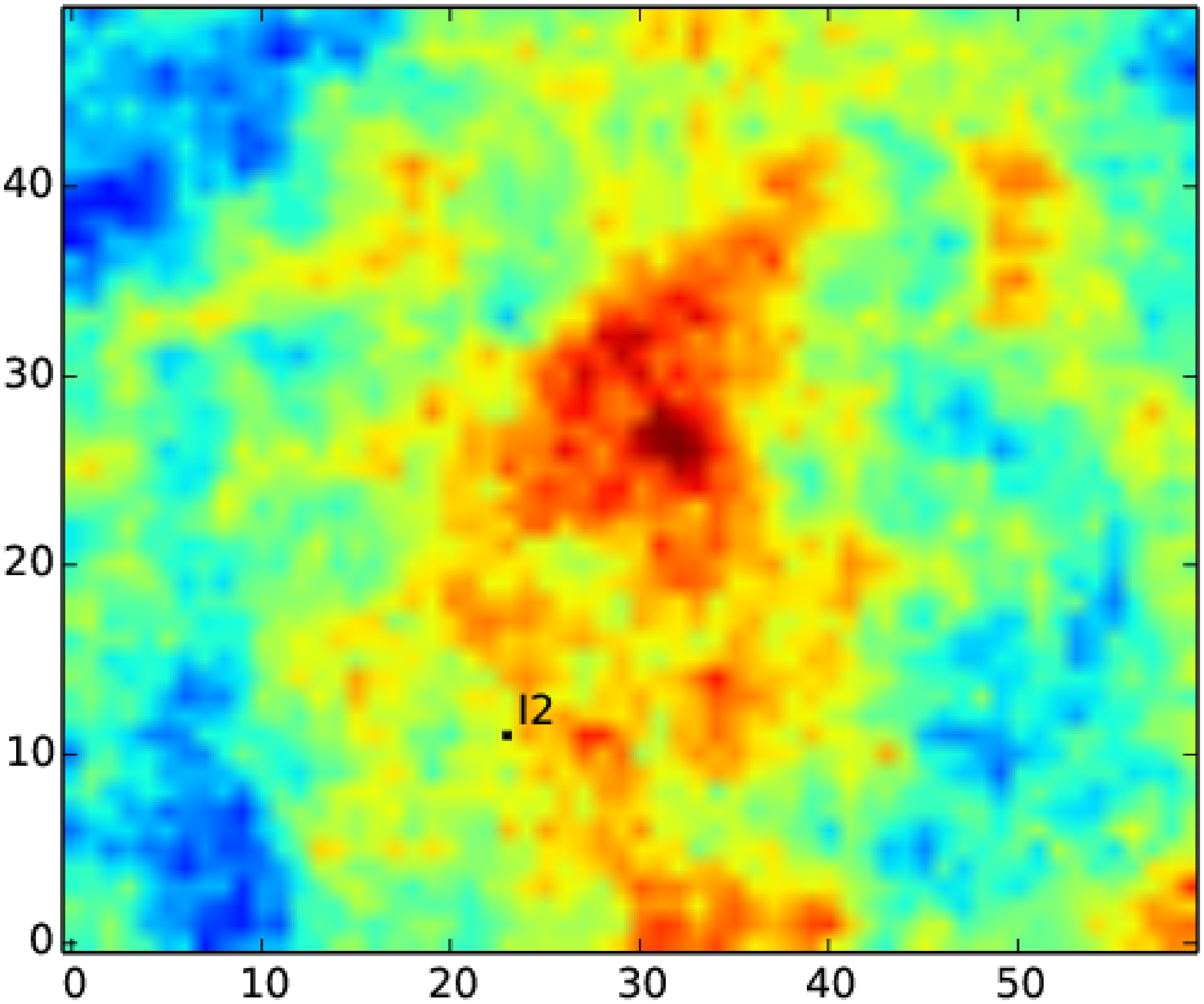} \\
\end{tabular}

\begin{tabular}{C{0.22\linewidth}C{0.22\linewidth}C{0.24\linewidth}C{0.24\linewidth}}
	\includegraphics[width=\linewidth]{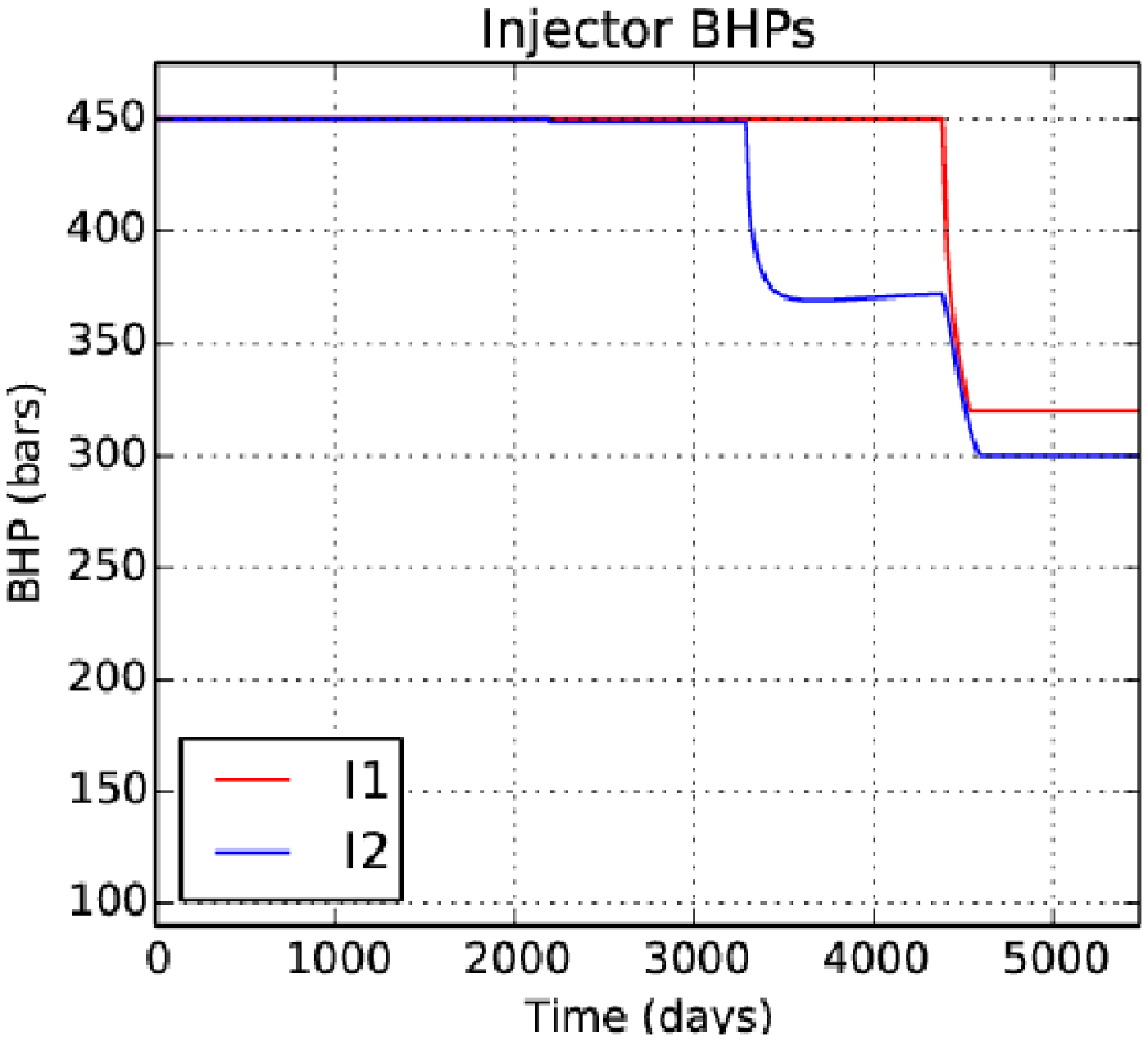}
&	\includegraphics[width=\linewidth]{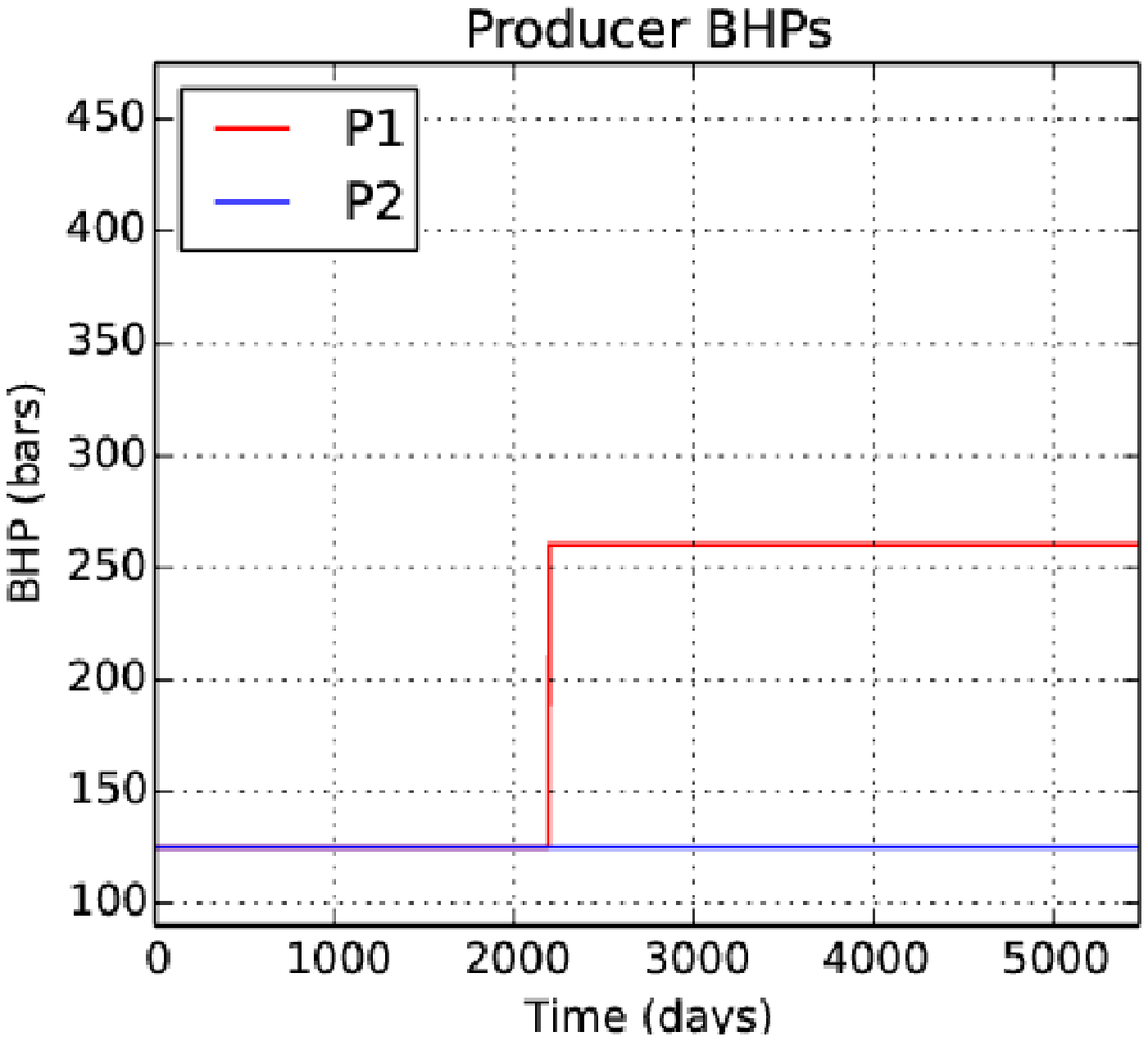}
&	\includegraphics[width=\linewidth]{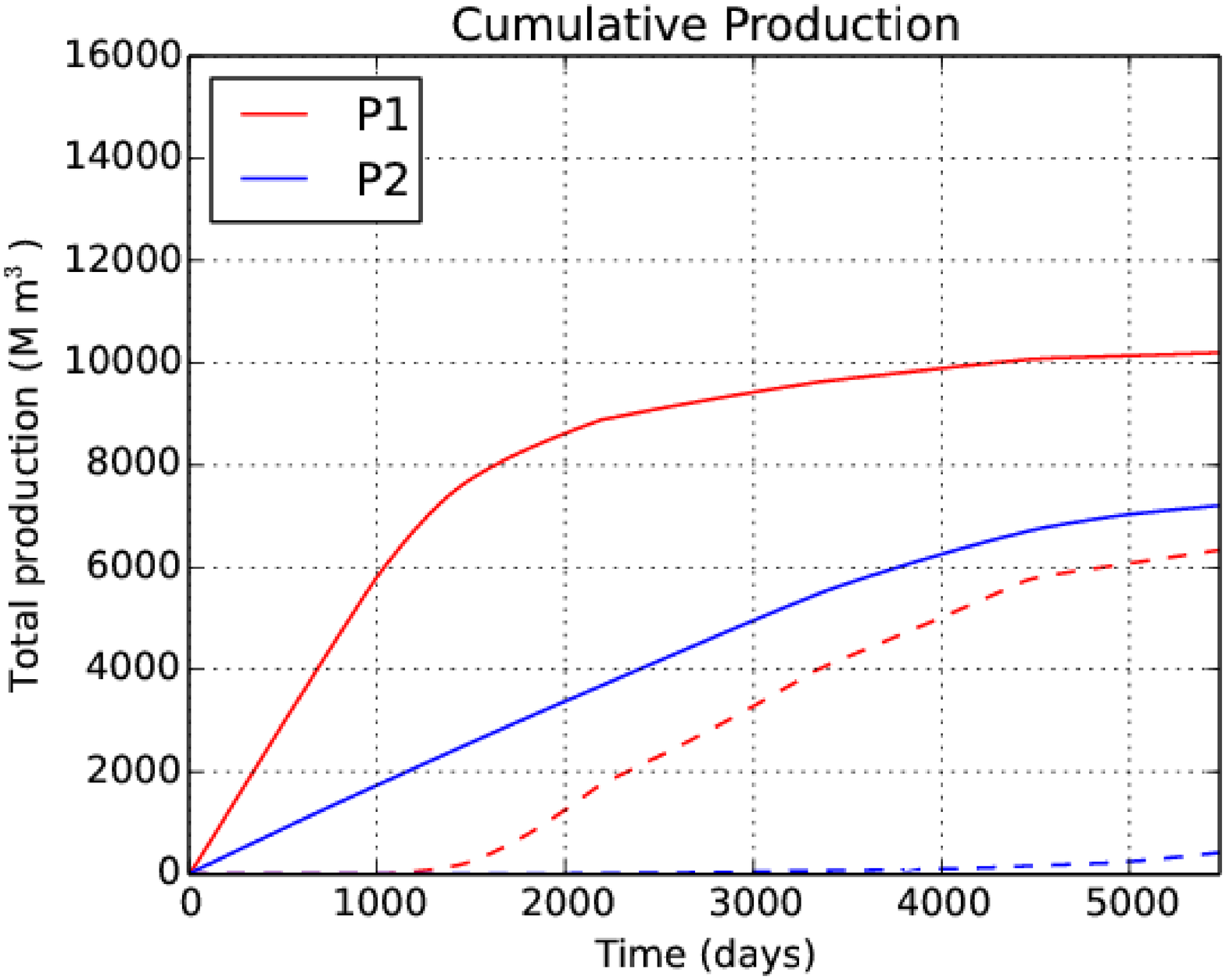}
&	\includegraphics[width=\linewidth]{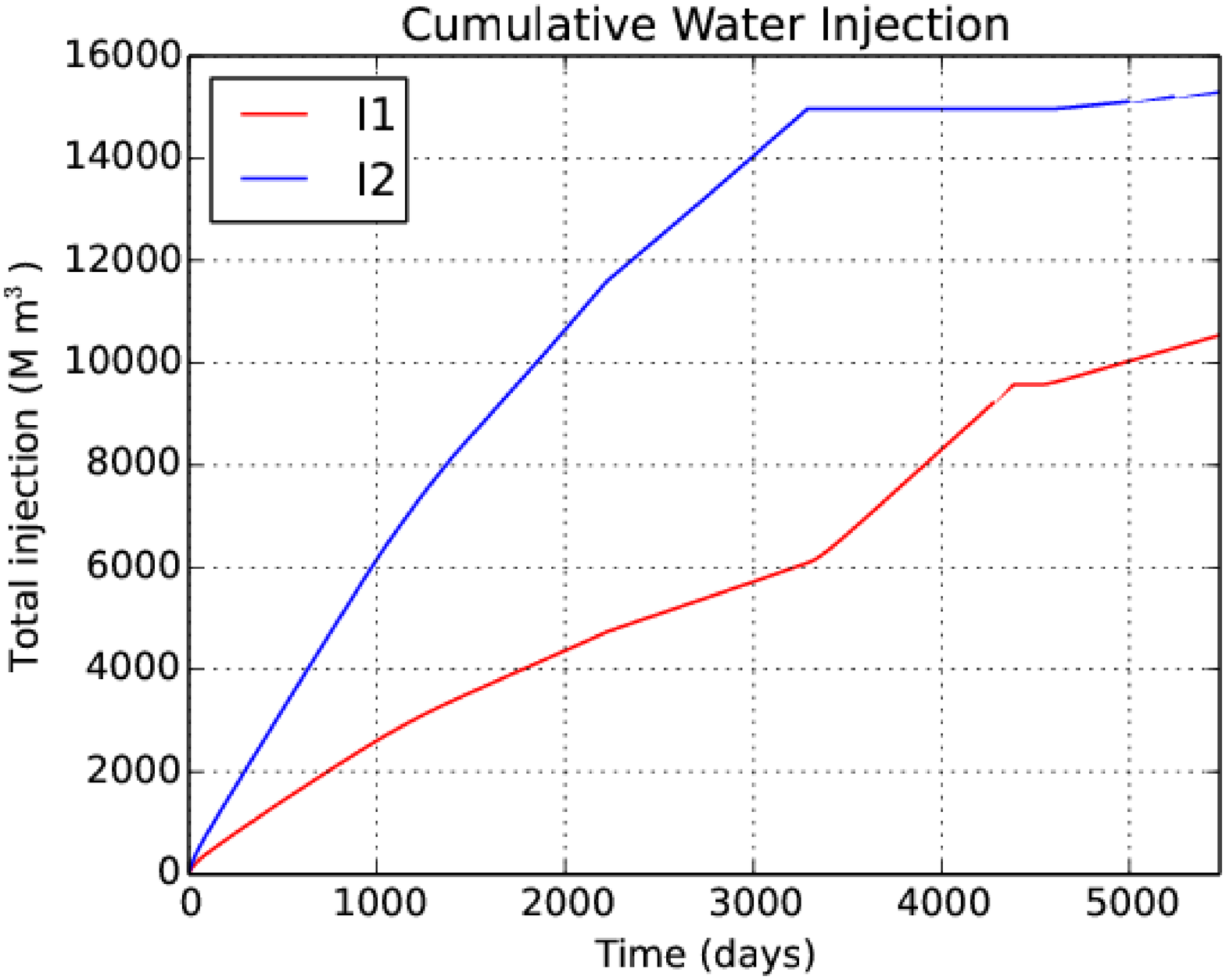} \\
\end{tabular}
\caption{Best solution found by Sequential-I for Case 3A, with NPV of \$4.850$\times 10^{9}$. Top row: permeability field (log scale) with well locations indicated. Bottom row: BHP controls for each injector and producer (left and centre-left plots), cumulative fluid production for each producer (centre-right) and cumulative water injection for each injector (right). The solid lines in the cumulative production plot indicate oil production; dashed lines indicate water production.}\label{F:Case3A_SEQ_best}
\end{figure*}

\begin{figure*}
\begin{tabular}{C{0.32\linewidth}C{0.32\linewidth}C{0.32\linewidth}}
{\bf Top Layer} &{\bf Middle Layer} &{\bf Bottom Layer} \\
\includegraphics[width=\linewidth]{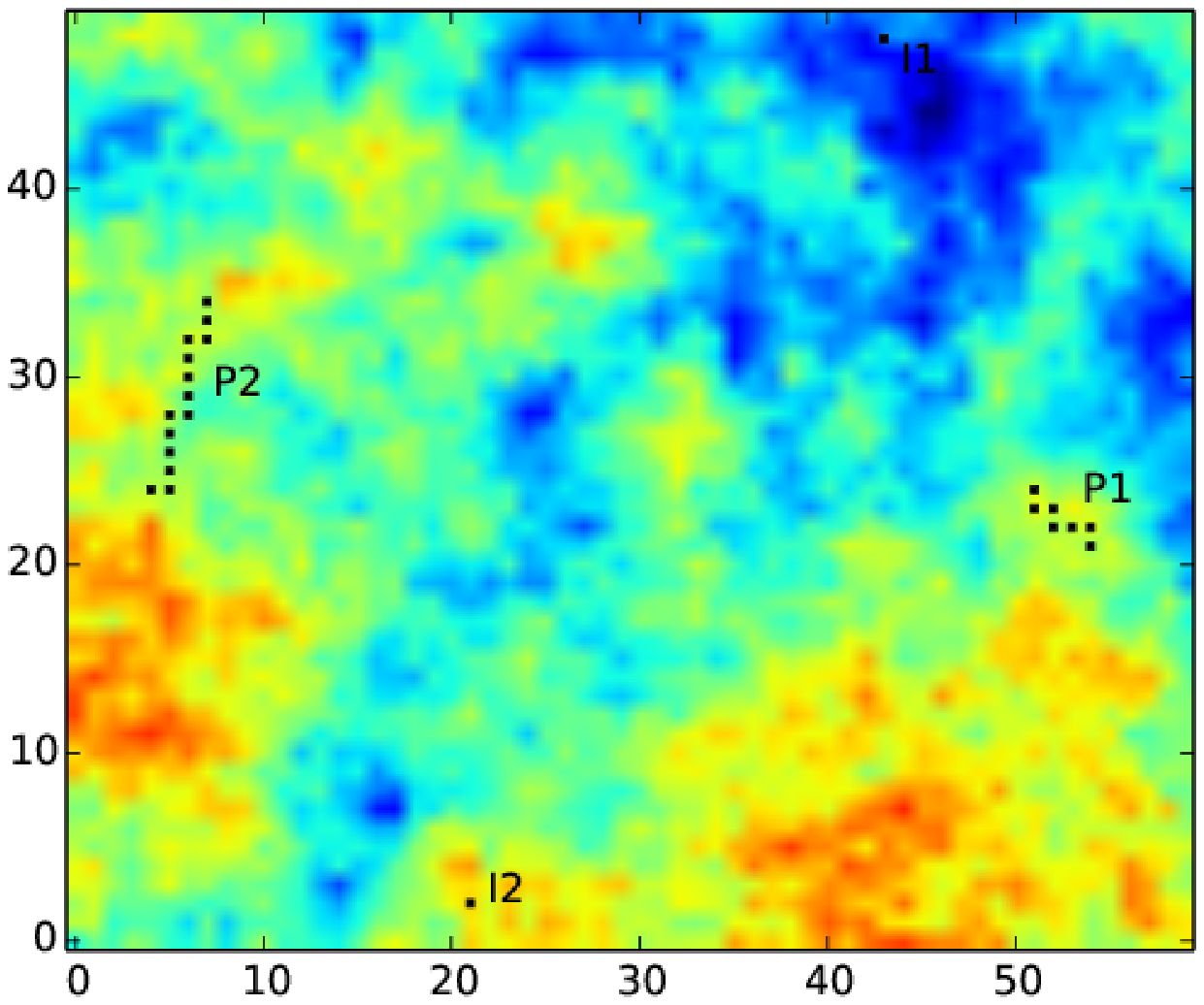}
&\includegraphics[width=\linewidth]{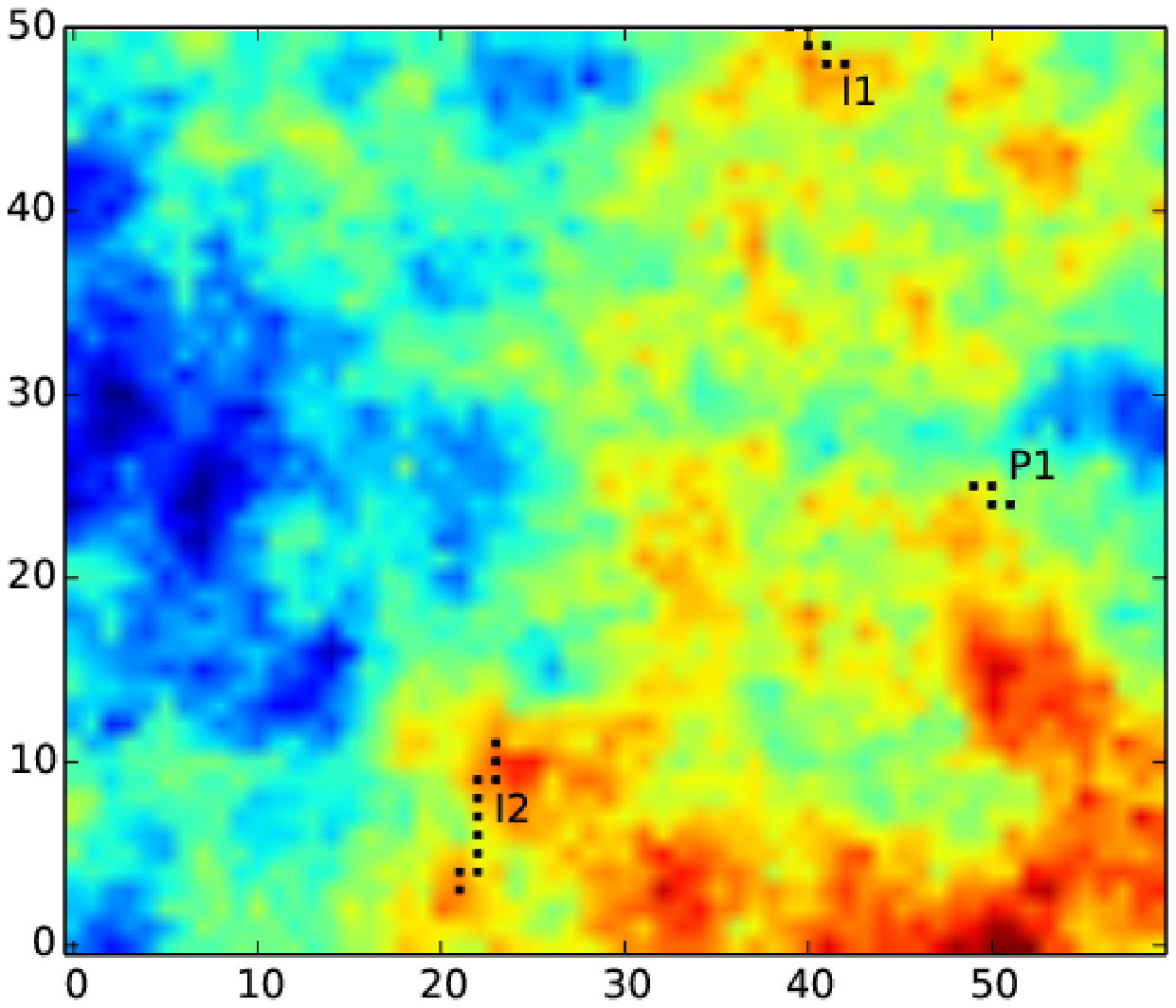}
&\includegraphics[width=\linewidth]{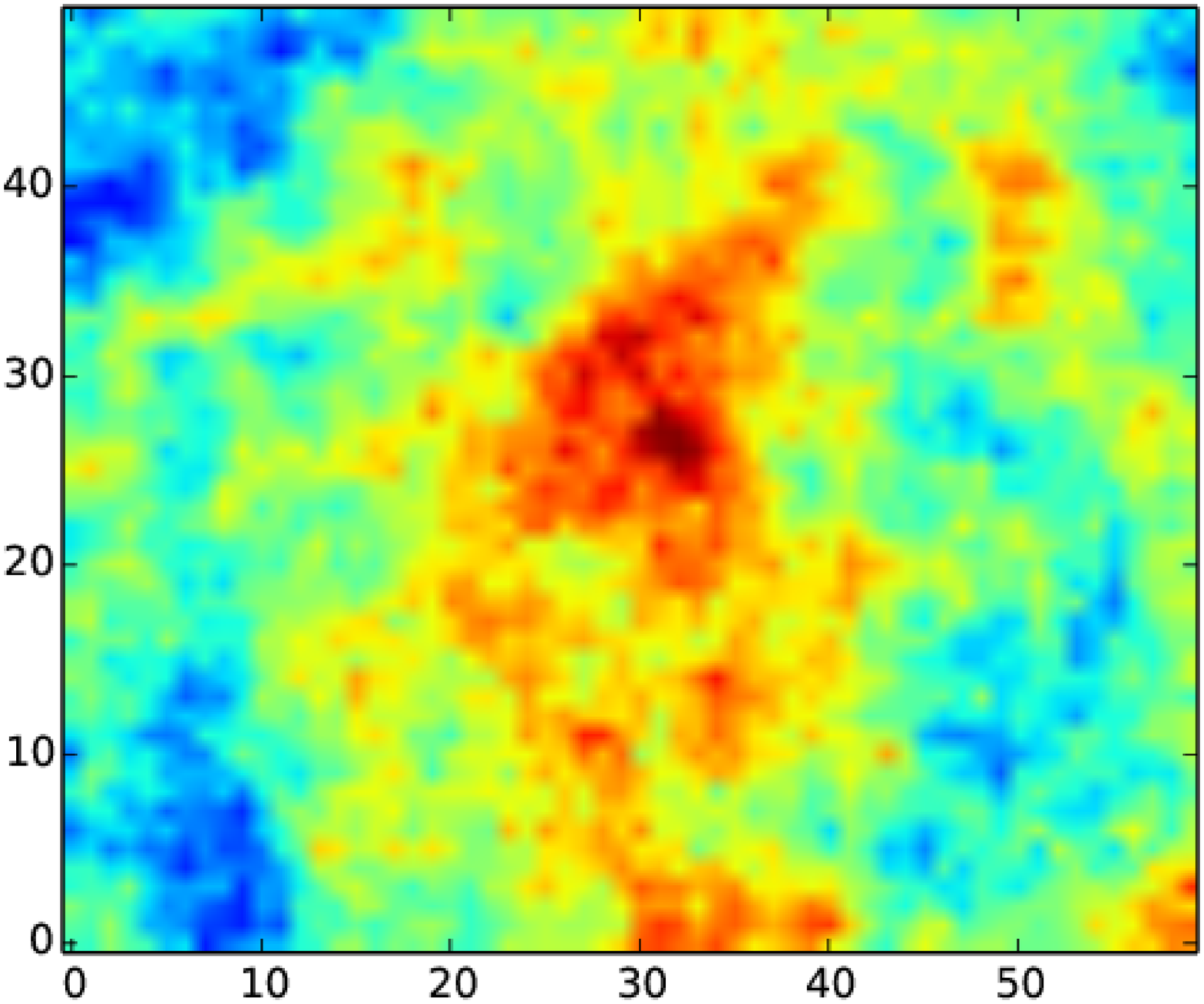} \\
\end{tabular}

\begin{tabular}{C{0.22\linewidth}C{0.22\linewidth}C{0.24\linewidth}C{0.24\linewidth}}
	\includegraphics[width=\linewidth]{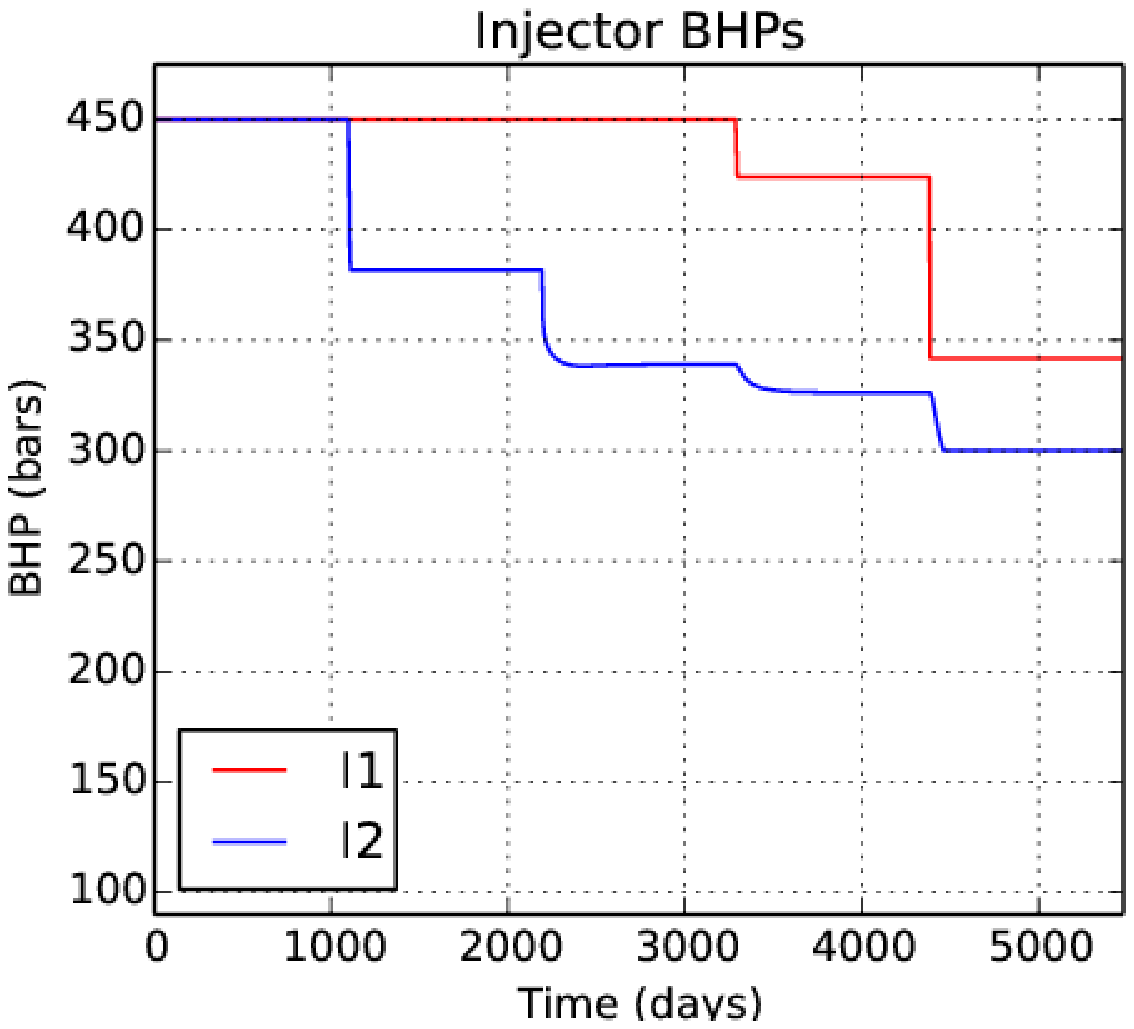}
&	\includegraphics[width=\linewidth]{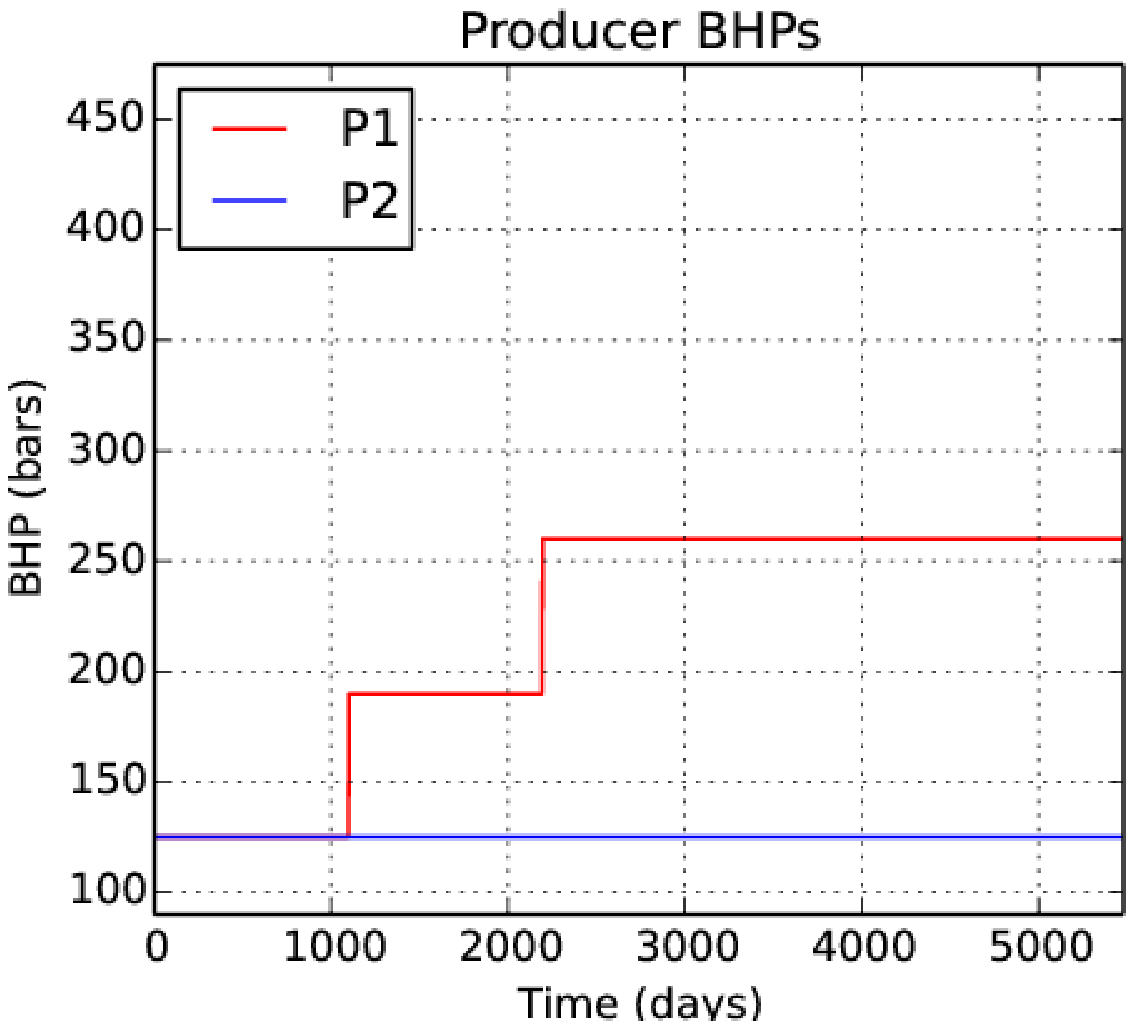}
&	\includegraphics[width=\linewidth]{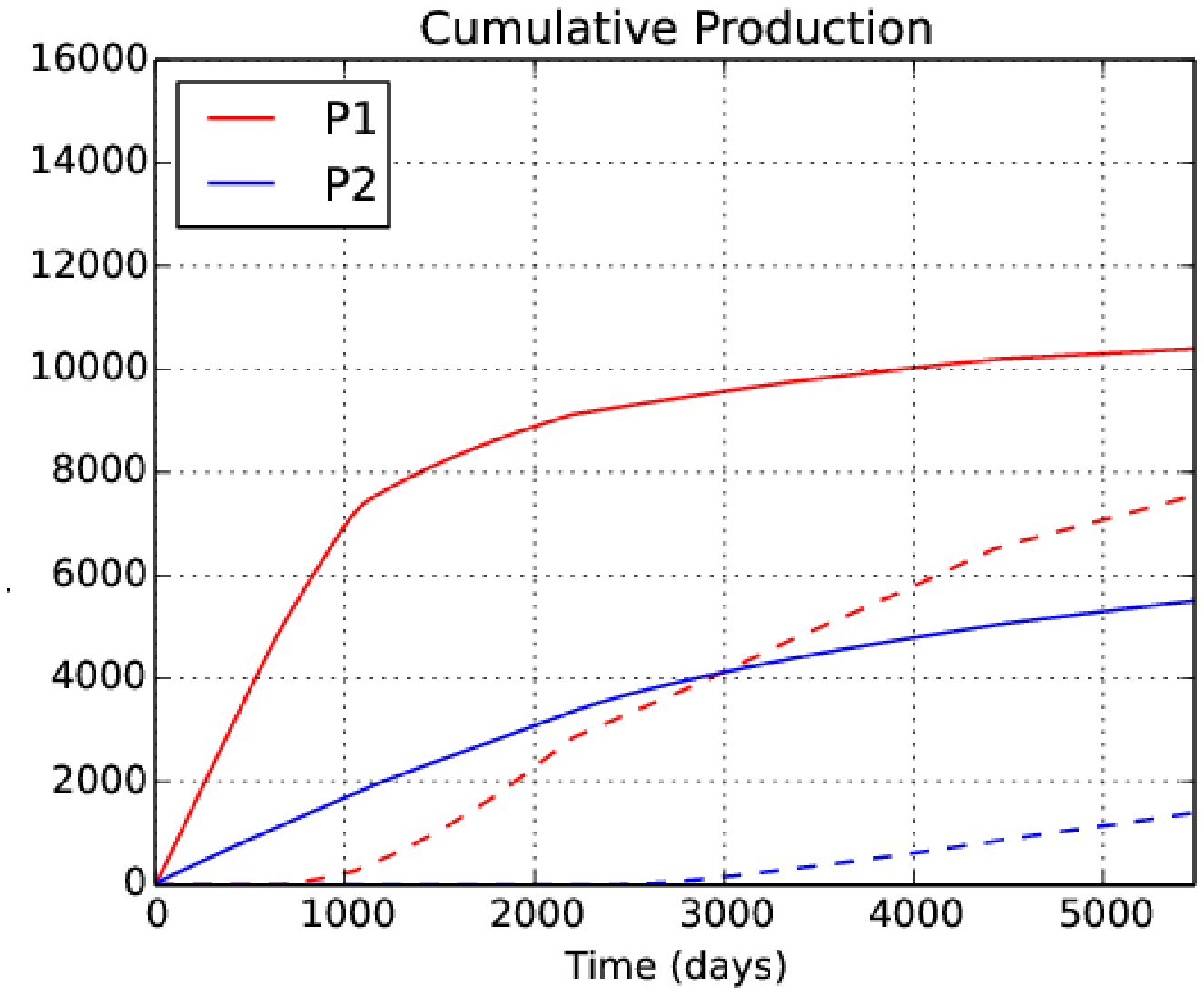}
&	\includegraphics[width=\linewidth]{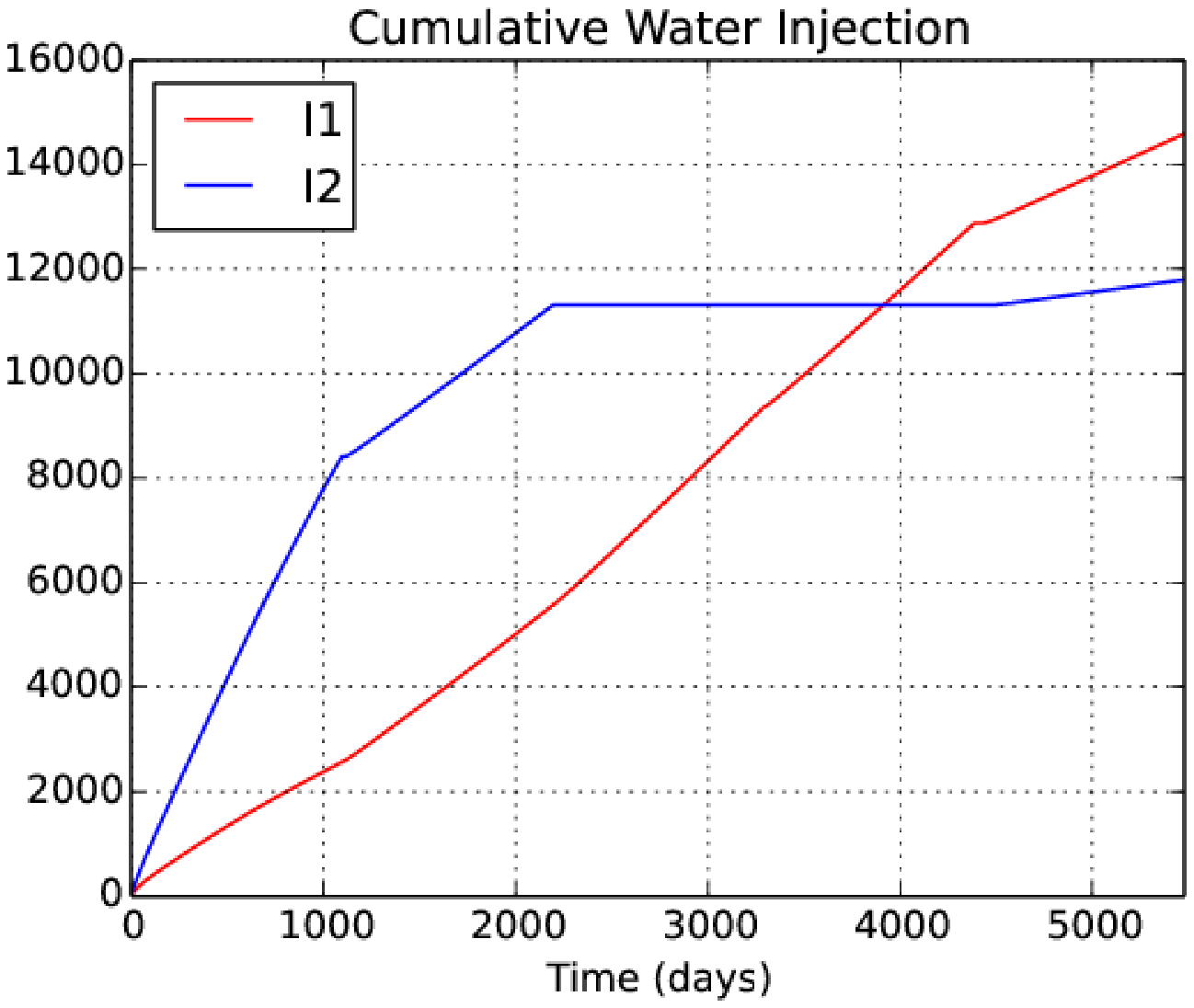} \\
\end{tabular}
\caption{Best solution found by MADS-PSO for Case 3A, with NPV of \$4.542$\times 10^{9}$. Plot and symbol meanings are the same as in \fref{F:Case3A_SEQ_best}.}\label{F:Case3A_MADS_PSO_best}
\end{figure*}

To gain further insight, we examine the best solutions found by the Sequential-I and simultaneous MADS-PSO approaches for Case 3A, which we will denote as Solution~1 and Solution~2, respectively. \fref{F:Case3A_SEQ_best} shows Solution~1 (NPV of \$4.850$\times10^{9}$), which was 6.8\% better than Solution~2, shown in \fref{F:Case3A_MADS_PSO_best}. The top row of both figures shows the $x$-$y$ permeability field for each layer of the reservoir, with cells perforated by each well indicated by black boxes. The bottom row of each figure shows the optimized BHP controls for each injector and producer, as well as the cumulative fluid production and cumulative water injection curves. We see from the top row of both figures that injectors and producers were placed in roughly the same area of the reservoir in both solutions, with the location and orientation of the injector labeled I2 being virtually the same in both solutions. The numerical results of the simulation, however, indicate that roughly 9.4\% more oil was produced in Solution~1 than in Solution 2, as well as nearly 25\% less water, which accounts for the discrepancy in NPV. This is also observable from the cumulative production curves shown in \fref{F:Case3A_SEQ_best} and \fref{F:Case3A_MADS_PSO_best}. The amount of water injected was roughly the same in both simulations. The  primary difference is that the injector labeled I1 and the producer labeled P2 in both figures are not placed as optimally in Solution~2. In particular, in Solution~1, I1 and P2 are oriented such that I1 drives oil towards P2, while in Solution~2, most of the injected water from I1 drives oil to P1. This results in less oil production at P2 and a greater amount of water production at P1, thereby reducing the NPV. We note that the control scheme determined in Solution~2 attempts to mitigate this effect by raising the BHP of P1 after only 3 years to reduce the flowrate. Since the wells are not optimally positioned, however, the control scheme only compensates to a limited extent. %{\bf longer wells, I1 and I2 adjusted}

%{\bf NB: are these results due to larger reservoir size meaning that max/min controls are a nearly optimal choice? should run a set of experiments with smaller reservoir size to see if we observe the same discrepancy since this might require varying controls more.}

\section{Conclusions} \label{S:concl}

In this paper we have considered two approaches to the joint well placement and control optimization problem. The first approach optimizes over all well placement parameters and well control parameters simultaneously, while the second approach consists of a first step that determines optimal well positionings subject to a simple control strategy, and a second step that uses the first solution as a starting point to optimize over all parameters. In our numerical experiments we gradually increase the number of parameters required to describe each well's position, while holding the number of control parameters fixed. Thus, in Experiments 1, 2 and 3, the positional parameters make up 28\%, 44\%, and 55\% of the total number of variables under consideration, respectively. We observe that in Experiment 1, the simultaneous approach provided significantly better results in one of two test cases; in Experiment 2, one of the two sequential approaches considered was able to provide comparable results to the simultaneous approach, and in Experiment 3, both sequential approaches significantly outperformed the fully simultaneous approach. 

Although these results are empirical, they suggest that as the well placement component of the joint optimization problem becomes more challenging, the sequential approach benefits from devoting more attention to that subproblem. While the fully simultaneous approach is, in principle, capable of finding any of the optimal solutions found by the sequential approaches, it appears that the inclusion of control parameters along with positional parameters in a single stage of optimization makes it more difficult to find these solutions. This finding may be broadly applicable to other applied optimization problems as well. If one is aware that certain variables are more crucial to the quality of a solution than others (as is the case with positional parameters in this problem), then there may be an advantage in devoting more attention to those variables, by splitting the optimization routine into stages. One caveat is that the choice of values for the other parameters which are held fixed during any stage of the optimization may be important. For example, in Experiment 2 of this paper we observed that the choice of fixed control scheme for the sequential approach (Sequential-I versus Sequential-II) had a significant impact on the effectiveness of the approach.

There are many possible avenues for future work. With respect to the optimization algorithms used, it is clear from the work presented that combining a stochastic algorithm such as PSO with the deterministic search provided by MADS is an effective strategy. Some potential areas of investigation are the use of other stochastic approaches such as differential evolution or covariance matrix adaptation to see if they offer any improvement in performance, as well as improved methods for handling nonlinear constraints. With respect to the production optimization problem, one can consider including additional decision parameters such as the number and type of well to drill, in addition to the scheduling of drilling operations~\citep[e.g.][]{IED14,ID14}. It is likely that introducing these variables will necessitate further development of the optimization approach, since some of them are categorical in nature. Along the lines of the work presented here, it will be particularly interesting to see whether dividing the optimization into several stages is still an effective strategy.  
\section*{Acknowledgements}
The authors acknowledge funding from the Natural Sciences and Engineering Research Council of Canada (NSERC), the Atlantic Canada Opportunities Agency (ACOA) and Research \& Development Corporation of Newfoundland and Labrador (RDC).

\end{document}